\documentclass[draft]{interact}

\usepackage{epstopdf}
\usepackage[caption=false]{subfig}

\usepackage[numbers,sort&compress]{natbib}
\usepackage{mathtools}
\usepackage[nameinlink, noabbrev]{cleveref}
\usepackage[shortlabels]{enumitem}
\bibpunct[, ]{[}{]}{,}{n}{,}{,}

\DeclarePairedDelimiterXPP\ind[1]{\mathbb{I}}( ){}{	#1}
\DeclarePairedDelimiterXPP\pk[1]{\mathbb{P}}\{ \}{}{ #1}
\DeclarePairedDelimiterXPP\E[1]{\mathbb{E}}\{ \}{}{	#1}

\newtheorem{theo}{Theorem}[section]
\newtheorem{lem}[theo]{Lemma}
\newtheorem{de}[theo]{Definition}
\def\EH#1{\textcolor{c40}{#1}}
\def\EH#1{#1}
\def\vv{v }

\def\bqn#1{ \begin{eqnarray} #1 \end{eqnarray}}
\def\MB{\mathcal{B}}
\def\KK{E}
\def\IF{\infty}
\def\LT{\left}
\def\RT{\right}
\def\rw{\rightarrow}
\def\N{\mathbb{N}}
\newcommand{\COM}[1]{}
\newcommand{\inr}{\in \R}
\newcommand{\ABs}[1]{ \left| #1 \right|}
\newcommand{\norm}[1]{\lVert  #1 \rVert }
\newcommand{\vk}[1]{{#1}}
\newcommand{\limit}[1]{\lim_{#1 \to   \infty}}
\newcommand{\BT}{\begin{theo}}
\newcommand{\ET}{\end{theo}}
\newcommand{\BEX}{\begin{example}}
\newcommand{\EEX}{\end{example}}
\newcommand{\abs}[1]{\left\lvert #1 \right\rvert}
\newcommand{\BD}{\begin{de}}
\newcommand{\ED}{\end{de}}
\newcommand{\BQN}{\begin{eqnarray}}
\newcommand{\EQN}{\end{eqnarray}}
\newcommand{\BQNY}{\begin{eqnarray*}}
\newcommand{\EQNY}{\end{eqnarray*}}
\newcommand{\R}{\mathbb{R}}
\newcommand{\BEL}{\begin{lem}}
\newcommand{\EEL}{\end{lem}}

\theoremstyle{plain}
\newtheorem{theorem}{Theorem}[section]

\newtheorem{corollary}[theorem]{Corollary}

\theoremstyle{definition}

\newtheorem{example}[theorem]{Example}

\theoremstyle{remark}

\begin{document}
	
	\title{Sojourns of locally self-similar Gaussian processes  }
	
	\author{
	\name{
	Svyatoslav ~M. Novikov\thanks{
		CONTACT Svyatoslav ~M. Novikov. Email:
		Svyatoslav.Novikov@unil.ch
	}
	}
	\affil{
	Department of Actuarial Science, University of Lausanne, Unil-Dorigny,
	1015 Lausanne, Switzerland}
	}
	
	\maketitle
	
	\begin{abstract}
		Given a Gaussian risk process $R(t)=u+c(t)-X(t),t\ge 0$,  the cumulative Parisian ruin probability on a finite time interval $[0,T]$ with respect to $L \geq 0$ is defined as the probability that the sojourn time that the risk process $R$ spends under the level 0 on this time interval $[0,T]$ exceeds $L$. In this contribution we derive exact asymptotic approximations of the cumulative Parisian ruin probability for a general class of Gaussian processes introduced in \cite{MR4206416} assuming that $X$ is locally self-similar. We illustrate our findings with several examples. As a byproduct we show that Berman's constants can be defined alternatively by a self-similar Gaussian process which could be quite different to the fractional Brownian motion.
	\end{abstract}
	

	\begin{amscode}
	Primary 60G15;
	Secondary 60G70
	\end{amscode}

	\section{Introduction}
	The Gaussian risk model $R(t)=u+ c(t)- X(t), t\ge 0$ is a classical object in modeling the risk of a certain insurance portfolio, see e.g., \cite{HP99,DJR20,debicki2020extremes} . Typically,  $X$ is a Brownian motion modelling the total claim amount, $u$ stands for the initial capital and $c(t)$ models the cumulative premium income up to time $t\ge 0$. The classical ruin time is the first time when the risk process $R$ drops below zero in a given interval, say $[0,T]$. The cumulative Parisian ruin, also known as the sojourn ruin is recently considered in the literature, see e.g., \cite{jasnovidov2020approximation,kriukov2022parisian,krystecki2022parisian,MR3457055} and the references therein. By definition, the cumulative Parisian ruin in the time interval $[0,T]$ occurs when the total sojourn time (occupation time) that $R$ spend below zero in the time interval $[0,T]$ exceeds a certain predefined threshold, say $L\ge 0$. More specifically, since the sojourn time that $R$ spends below 0 is given by 
	$$\int_{0}^T \ind{R(t) < 0}dt,   $$
	with $\ind{\cdot}$ the indicator function,  the sojourn probability is thus defined by 
	$$ p_X(T,L)= \pk*{\int_{0}^T \ind{X(t) - c(t) > u}dt > L  }.$$
	Clearly, when $L=0$, we retrieve the classical ruin probability on the finite time horizon $[0,T]$.%
	\\
	The asymptotics of sojourn times has been the main topic of many research papers by Berman, see e.g., \cite{MR803245, Berman92} where the stationary and locally stationary cases are discussed. Recently, \cite{KEZX17,debicki2023sojourn} 
	derived results for the tail asymptotics of the sojourn times of general Gaussian processes and fields. 
	
	A typical assumption in the literature on the  Gaussian process $X$ is that 
	the variance function $\sigma^2(t),t\in [0,T]$  attains its maximum at some unique point, say $t_0=0$ and further 
	\bqn{ 
		\label{sigmaE}
		\sigma(t_0)-\sigma(t) \sim bt^{\beta} , \quad t\downarrow   0,
	}
	with some $b, \beta$ positive constants; here $\sim$ stands for asymptotic equivalence when the argument tends to 0 or $\infty$ (depending on the context).\\
	The assumption on the variance function is accompanied by an asymptotic assumption on the correlation function \cite{Pit96}.
	Specifically, it is an asymptotic condition on the correlation $r_X$ of $X$ at zero, saying that 
	\bqn{ 
		\label{PickE}
		1-r_{X}(t,s) \sim a V_Y(t,s)  , \quad s,t\downarrow  0
	}
	for some $a>0, \kappa \in (0,2]$, where 
	$$V_Y(t,s)=\mathrm{Var}(Y(t)-Y(s))=\abs{t-s}^\kappa, \quad s,t\ge 0$$ is the variogram of a
	centered fractional Brownian motion (fBm).

	In this contribution, we shall consider the locally self-similar risk model, i.e., $X$ is locally self-similar as introduced in \cite{tabis,MR4206416}. This definition allows for a more general $Y(t),t\ge 0$, which is 
	again a centered Gaussian process with a.s. continuous sample paths satisfying (set below $V_Y(t)=V_Y(t,1)=\mathrm{Var}(Y(1)-Y(t))$)
	\begin{enumerate}[S1)]
		\item \label{LS1}
		$Y$ is self-similar with index $\alpha / 2>0$ and $\sigma_Y^2(1)=1$;
		
		\item \label{LS2} there exist $\kappa \in(0,2]$ and $c_Y>0$ such that $V_Y(1-h)\sim c_Y|h|^\kappa$ as $h \rightarrow 0$.
	\end{enumerate}
	For notational symplicity, $Y \in \mathbf{S}\left(\alpha, \kappa, c_Y\right)$ shall stand as an abbreviation that $Y$ satisfies \Cref{LS1}-\Cref{LS2}. \\
	\BD A centered Gaussian process $X(t), t \geq 0$ is said to be locally self-similar at $0$ with parameters $a$ and $Y$, if \eqref{PickE} holds for some $a>0, Y \in \mathbf{S}\left(\alpha, \kappa, c_Y\right)$. 
	\ED

	Brief organisation of the paper: In our first result we show that the Berman constant defined with respect to a standard Brownian motion, can be
	alternatively defined with respect to a self-similar process, see Thm 2.2. Then we obtain the tail asymptotics of the cumulative Parisian ruin (sojourn time) for locally self-similar Gaussian processes. Our findings are illustrated by several examples. 
	Auxiliary results are presented in Section 4 followed by the proofs postponed to Section 5.

	\section{Main Results}
	In this section we shall investigate the cumulative Parisian ruin considering as trend function $c(t)= d t^{\gamma}$ for some $d,\gamma$ positive constants. More general trends can be dealt with similar techniques and will therefore not be considered here. More specifically, we shall investigate the asymptotics of  
	$$ p_X(T,L_u)=\pk*{\int\limits_{0}^{T} \mathbb{I}\left(X\left(t\right)-dt^{\gamma}>u\right) dt > L_u} $$
	as the initial capital $u \to \IF$ for constant $L_u$ dependent  on $u$ and  $X$ being a locally self-similar Gaussian process.

	In \cite[Thm 3.4]{MR4206416} under some conditions a relation between the Pickands constant of a self-similar process $Y$ and the Pickands constant of a fractional Brownian motion $B_{\kappa}$ was proved.

	Recently the properties of Berman constants were investigated in \cite{dkebicki2022berman}. As therein, the Berman  constant with respect to a Gaussian process $\zeta$ can be defined as  
	\BQNY
	\MB_{\zeta}^{h}\left( x\right)= \lim_{T\to\infty} \frac{1}{T} \int_\R  \pk*{ \int_{ 0 }^T \mathbb{I}\left(\sqrt{2}\zeta(t)-\mathrm{Var}\left(\zeta(t)\right)-h(t) + y >0\right)  dt >x } e^{-y} dy.
	\EQNY
	For convenience we will denote $\MB_{\zeta}^{0}(x)$ by $\MB_{\zeta}(x)$. We also introduce
	$$\MB_{\zeta}^{h}(x,E) =
	\int_\R  \pk*{ \int_{ E } \mathbb{I}\left(\sqrt{2}\zeta(t)-\mathrm{Var}\left(\zeta(t)\right)-h(t) + y >0\right)  dt >x } e^{-y} dy.$$
	We show below that as in \cite[Thm 3.4]{MR4206416} the Berman constants can be defined directly with respect to a standard fractional Brownian motion, which is a remarkable result.\\ 
	There is an additional technical difficulty, compared to \cite[Thm 3.4]{MR4206416} because we cannot apply the Slepian inequality.
	
	\BT\label{last}
	If $\widehat{Y} \in \mathbf{S}\left(\kappa,\kappa,c_{\widehat{Y}}\right)$ and $x \geq 0$, then \BQN \MB_{\widehat{Y}}\left(x\right)=c_{\widehat{Y}}^{1/\kappa}\MB_{B_{\kappa}}\left(c_{\widehat{Y}}^{1/\kappa}x\right)\in (0,\infty).
	\EQN
	\ET
	Remark that \Cref{last} does not extend to $\MB_Y(x)$ with $Y \in \mathbf{S}(\alpha,\kappa,c_Y)$ in the same way as \cite[Thm 3.4]{MR4206416}, since \cite[Lm 4.1 (ii)]{MR4206416} cannot be generalized to $\MB_Y(x)$. Indeed, setting $\widehat{Y}(t) = Y(t^{\kappa/\alpha})$ as in \cite{MR4206416}, one obtains 
	\BQNY
	&&
	\MB_Y(x)\\ && = 
	\lim_{T\to\infty} \frac{1}{T} \int_\R  \pk*{ \int_{ 0 }^T \mathbb{I}\left(\sqrt{2}\widehat{Y}(t^{\alpha/\kappa})-\mathrm{Var}\left(\widehat{Y}(t^{\alpha/\kappa})\right) + y >0\right)  dt >x } e^{-y} dy
	\\
	&& =
	\lim_{T\to\infty} \frac{1}{T} \int_\R  \pk*{ \int_{ 0 }^{T^{\alpha/\kappa}} \mathbb{I}\left(\sqrt{2}\widehat{Y}(t)-\mathrm{Var}\left(\widehat{Y}(t)\right) + y >0\right) \frac{\kappa}{\alpha}t^{\kappa/\alpha-1} dt >x } e^{-y} dy,
	\EQNY
	which is not expressible in terms of $\MB_{\widehat{Y}}(\cdot)$.
	
	
	We extend \cite[Thm 4.2]{tabis} to the case of sojourn functional. Set in the following 
	$$\hat{\beta}=\frac{\beta \kappa}{\alpha},\;\hat{\gamma}=\frac{\gamma \kappa}{\alpha}.
	$$ Write next $\Psi$ for the survival function of a standard normal random variable.
	
	\BT\label{mainth}
	If $X(t), t \geq 0$ is locally self-similar at $0$ with parameters $a$ and $Y$ which satisfies \eqref{sigmaE} then we have 
	\BQN
	\label{mainasymp}
	\pk*{\int\limits_{0}^{T} \mathbb{I}\left(X\left(t\right)-dt^{\gamma}>u\right) dt > L_u}\sim c u^{p}\Psi\left(u\right), \quad u\to \IF,
	\EQN
	where 
	
	\begin{enumerate} [(i)]
		\item 
		If $\alpha < \min\left(\beta,2\gamma\right),\;\alpha \leq \kappa$, then for $L_u=Lu^{-2/\alpha+\left(\frac{\kappa-\alpha}{\kappa}\right)\left(\frac{2}{\alpha}-\max\left(1/\gamma,2/\beta\right)\right)}$, $L\geq 0$ it holds that \\$p=2/\kappa-\max\left(2/\hat{\beta},1/\hat{\gamma}\right)$,
		
		\begin{multline*}
		c=\left(ac_{Y}\right)^{1/\kappa}
		\int\limits_{0}^{\infty} 
		\exp\left(-d\mathbb{I}
		\left(2\gamma \leq \beta\right) z^{\gamma}
		-b\mathbb{I}\left(\beta \leq 2\gamma\right) z^{\beta}\right)
		\\
		\times
		z^{\alpha/\kappa-1}\MB_{B_{\kappa}}\left(L\left(ac_{Y}\right)^{1/\kappa} z^{\alpha/\kappa-1}\right) dz.
		\end{multline*}
		
		\item  If $\alpha<\min\left(\beta,2\gamma\right),\;\alpha>\kappa$, then for $L_u=Lu^{-2/\alpha-\epsilon}$,
		\\ $\epsilon \in (0,\frac{\alpha-\kappa}{\alpha}\left(2/\kappa-\max\left(2/\hat{\beta},1/\hat{\gamma}\right)\right)]$, $L>0$ it holds that $p=\frac{\epsilon\alpha}{\alpha-\kappa}$,
		
		\begin{multline*}
		c= \left(ac_{Y}\right)^{1/\kappa}
		\int\limits_{0}^{\infty} 
		\exp\Bigg(-d\mathbb{I}
		\left(\epsilon=\frac{\alpha-\kappa}{\alpha}\left(2/\kappa-1/\hat{\gamma}\right)\right) z^{\gamma}
		\\
		-b\mathbb{I}\left(\epsilon=\frac{\alpha-\kappa}{\alpha}\left(2/\kappa-2/\hat{\beta}\right)\right) z^{\beta}\Bigg)
		z^{\alpha/\kappa-1}\MB_{B_{\kappa}}\left(L\left(ac_{Y}\right)^{1/\kappa} z^{\alpha/\kappa-1}\right) dz.
		\end{multline*}
		If further $\epsilon =\frac{\alpha-\kappa}{\alpha}\left(2/\kappa-\max\left(2/\hat{\beta},1/\hat{\gamma}\right)\right)$, then the same holds with $L=0$.
		
		\item If $\alpha \geq \min\left(\beta,2\gamma\right)$ or $\alpha>\kappa$, then for $L_u=Lu^{\min\left(-2/\alpha,-2/\beta,-1/\gamma\right)}$, $L>0$ it holds that $p=0$, \\$$c=\MB_{Y\mathbb{I}\left(\alpha \leq \min\left(\beta,2\gamma\right)\right)}^{h}
		\left(L
		a^{1/\kappa},[0,\infty)\right)$$ with $$h\left(t\right)=a^{-\beta/\alpha}bt^{\beta}\mathbb{I}\left(\beta\leq \min\left(\alpha,2\gamma\right)\right)+
		a^{-\gamma/\alpha}
		dt^{\gamma}\mathbb{I}\left(2\gamma\leq \min\left(\alpha,\beta\right)\right).$$
		
		If $\alpha \geq \min\left(\beta,2\gamma\right)$, then the same holds for $L=0$.
		
		In particular, the corresponding Berman constants are positive and finite.
	\end{enumerate}
	\ET
	\begin{corollary}\label{corselfsim}
		If $Y \in \mathbf{S}\left(\alpha, \kappa, c_Y\right)$, $V_Y(x)$ is decreasing in some neighbourhood of $0$, $V_Y(x)$ attains its maximum on $[0,1]$ at the unique point $x=0$ and $1-V_Y\left(x\right)
		\sim R x^{\beta}$ with $\beta \geq 1$, $\beta>\alpha/2$, and in addition, the function $\frac{\frac{\partial}{\partial x}V_Y(x)}{x^{\beta-1}}$ is bounded on $(0,1]$,
		then
		the asymptotics from Theorem \ref{mainth} hold for $X\left(t\right)=Y\left(1\right)-Y\left(t\right)$ for $a=1/2$, $b = R/2$.
	\end{corollary}

	\section{Examples}
	In this section we consider different possible examples of $Y$ in Corollary \ref{corselfsim}.
	
	Remark that if $\alpha > \min\left(\beta,2\gamma\right)$, then
	\BQNY
	c&=&\MB_0^{h}\left(L a^{1/\kappa},[0,\infty)\right)=
	\exp\left(-h\left(L a^{1/\kappa}\right)\right)
	\\
	&=&a^{\beta/\kappa-\beta/\alpha} b t^{\beta} \mathbb{I} \left(\beta \leq \min(\alpha,2\gamma)\right)+
	a^{\gamma/\kappa-\gamma/\alpha} d t^{\gamma} \mathbb{I} \left(2 \gamma \leq \min(\alpha,\beta)\right).
	\EQNY
	
	Below
	$R_Y\left(t,s\right),t,s\in[0,T]$ stands for the covariance function of
	$Y\left(t\right), t \in [0,T]$.
	
	\BEX
	For $\alpha \in (1,2)$ consider the covariance function
	$$
	R_Y(t,s)=\frac{(t+s)^{\alpha}-|t-s|^{\alpha}}{2^{\alpha}}.
	$$
	
	In view of \cite{tabis} $\beta=1,R=\alpha 2^{2-\alpha}$, $Y \in \mathbf{S}(\alpha,\alpha,2^{1-\alpha})$. Let $\gamma=\beta/2=1/2$, then \eqref{mainasymp} holds with $L \geq 0$,
	$L_u=Lu^{-2}$, $p=0$, $c=e^{-\alpha 2^{1-\alpha}L-dL^{1/2}}$.
	\EEX
	
	\BEX
	Sub-fractional Brownian motion with parameter $\alpha \in (0,2)$ is a centered Gaussian process with covariance function
	$$
	R_Y(t,s)=\frac{1}{2-2^{\alpha-1}}
	(t^{\alpha}+s^{\alpha}-
	\frac{(t+s)^{\alpha}+|t-s|^{\alpha}}
	{2}).
	$$
	In view of \cite{tabis} $\beta=\alpha,R=\frac{2^{\alpha-1}}{2-2^{\alpha-1}},Y \in \mathbf{S}(\alpha,\alpha, (2-2^{\alpha-1})^{-1})$. Let $\gamma=\beta/2$, then \eqref{mainasymp} holds with $L \geq 0$,
	$L_u=Lu^{-2/\alpha}$, $p=0$ and $c=\MB_{Y}^{h}
	\left(2^{-1/\alpha}L,[0,\infty)\right)$, where $h(t)=\frac{2^{\alpha-1}}{2-2^{\alpha-1}}t^{\alpha}+\sqrt{2}dt^{\alpha/2}$.
	\EEX
	
	\BEX
	Negative sub-fractional Brownian motion $Y$ with parameter $\alpha \in(2,4]$ is a centered Gaussian process with covariance function
	$$
	R_Y(t, s)=\frac{1}{2^{\alpha-1}-2}\left(\frac{(t+s)^\alpha+|t-s|^\alpha}{2}-t^\alpha-s^\alpha\right) .
	$$
	In view of \cite{tabis}	
	$\beta=2,R=\frac{\alpha(\alpha-1)}{2^{\alpha-1}-2},Y\in \mathbf{S}(\alpha,2,\frac{\alpha(\alpha-1)2^{\alpha-3}}{2^{\alpha-1}-2})$. Let $\gamma=\beta/2=1$, then \eqref{mainasymp} holds with $L \geq 0$, $L_u=Lu^{-1}$, $p=0$ and
	$$c=\exp\left(-2^{2/\alpha}\frac{\alpha(\alpha+1)}{2^{\alpha+1}-8}L^2-2^{1/\alpha-1/2}dL\right).$$
	\EEX
	\BEX Weighted fBm $Y$ with parameters $\kappa \in(0,2]$ and $a > 1$ is a centered Gaussian process with covariance function
	$$
	R_Y(t, s)=\frac{\Gamma(a+\kappa)}{2 \Gamma(a) \Gamma(\kappa)} \int_0^{\min (t, s)} u^{a-1}\left((t-u)^{\kappa-1}+(s-u)^{\kappa-1}\right) d u .
	$$
	In view of \cite{tabis}
	$\beta = a,R=\frac{\Gamma(a+\kappa)}{\Gamma(a+1)\Gamma(\kappa)},Y\in \mathbf{S}(a+\kappa-1,\kappa,\frac{\Gamma(a+\kappa)}{\Gamma(a)\Gamma(\kappa+1)})$. Let $\gamma=\frac{\beta}{2}=\frac{a}{2},\alpha=a+\kappa-1.$
	
	If $\kappa>1$, then $\alpha>\beta$ and \eqref{mainasymp} holds with $L_u=Lu^{-2/a}$, $p=0$ and $$c=\exp\left(-2^{\frac{1-\kappa}{a+\kappa-1}-\frac{a}{\kappa}}\frac{\Gamma(a+\kappa)}{\Gamma(a+1)\Gamma(\kappa)}L^a-2^{\frac{a}{2(a+\kappa-1)}-\frac{a}{2\kappa}} dL^{a/2}\right).$$
	
	If $\kappa =1$, then $Y(t)$ has the same distribution as $B_1(t^a)$, and \eqref{mainasymp} holds with $L \geq 0$, $L_u=Lu^{-2/a}$, $p=0$ and $c=\MB_{Y}^{h}
	\left(2^{-1/\kappa}L,[0,\infty)\right)$, where $h(t)=t^a+\sqrt{2}dt^{a/2}$.
	
	If $\kappa<1$, then $\alpha<\beta$, $\alpha>\kappa$, hence, \eqref{mainasymp} holds with
	$$L_u=Lu^{-\frac{2}{a+\kappa-1}-\epsilon},
	\quad
	\epsilon \in \Big(0,\frac
	{2(a-1)(1-\kappa)}
	{\kappa a (a+\kappa-1)}\Big],\quad p=\epsilon \frac{a+\kappa-1}{a-1}$$ 
	and
	
	\begin{multline*}
		c=\left(\frac{\Gamma\left(a+\kappa\right)}{2\Gamma\left(a\right) \Gamma\left(\kappa+1\right)}\right)^{1/\kappa}
		\int\limits_{0}^{\infty} 
		\exp\left(-dz^{a/2}- 
		\frac{\Gamma\left(a+\kappa\right)}
		{2\Gamma\left(a+1\right)\Gamma\left(\kappa\right)} z^{a}\right) 
		\\
		\times
		z^{\frac{a-1}{\kappa}}
		\MB_{B_{\kappa}}\left(L\left(\frac{\Gamma\left(a+\kappa\right)}{2\Gamma\left(a\right) \Gamma\left(\kappa+1\right)}\right)^{1/\kappa}z^{\frac{a-1}{\kappa}}\right)dz,
	\end{multline*}
	with $L \geq 0$ if $\epsilon = \frac
	{2(a-1)(1-\kappa)}
	{\kappa a (a+\kappa-1)}$,
	
	$$
	c=\left(\frac{\Gamma\left(a+\kappa\right)}{2\Gamma\left(a\right) \Gamma\left(\kappa+1\right)}\right)^{1/\kappa}
	\int\limits_{0}^{\infty} 
	z^{\frac{a-1}{\kappa}}
	\MB_{B_{\kappa}}\left(L\left(\frac{\Gamma\left(a+\kappa\right)}{2\Gamma\left(a\right) \Gamma\left(\kappa+1\right)}\right)^{1/\kappa}z^{\frac{a-1}{\kappa}}\right)dz,
	$$
	with $L>0$ if $0 <\epsilon <\frac
	{2(a-1)(1-\kappa)}
	{\kappa a (a+\kappa-1)}$,
	
	$$
	c=\MB_{Y}^{h}
	\left(2^{-1/\kappa}L,[0,\infty)\right),
	$$
	where $h(t)= -2^{\frac{1-\kappa}{a+\kappa-1}}\frac{\Gamma(a+\kappa)}{\Gamma(a+1)\Gamma(\kappa)}t^a-2^{\frac{a}{2(a+\kappa-1)}} dt^{a/2}$
	with $L>0$ if $\epsilon = 0$.
	\EEX
	
	\BEX
	Integrated fBm $Y$ with parameter $\alpha \in(0,2)$ is a Gaussian process defined as
	$$
	Y(t)=\sqrt{\alpha+2} \int_0^t B_\alpha(s) d s .
	$$
	Its covariance function is of the form
	$$
	R_Y(t, s)=\frac{(\alpha+2)\left(s^{\alpha+1} t+s t^{\alpha+1}\right)+|t-s|^{\alpha+2}-t^{\alpha+2}-s^{\alpha+2}}{2(\alpha+1)} .
	$$
	In view of \cite{tabis} $Y \in \mathbf{S}(\alpha+2,2,\alpha+2)$ and
	$\beta=\alpha+1$ if $\alpha \leq 1$ and $\beta = 2$, if $\alpha \geq 1$. Let $\gamma = \beta/2$, then for some $b>0$ \eqref{mainasymp} holds with $L_u=Lu^{-\max\left(1,\frac{2}{\alpha+1}\right)}$, $p=0$ and $$c=\exp\left(-d \cdot
		2^{\min\left((\alpha+1)/2,1\right)\cdot\left(\frac{1}{\alpha+2}-\frac{1}{2}\right)}
		L^{\min\left((\alpha+1)/2,1\right)}-b
		L^{\min\left(\alpha+1,2\right)}\right).$$
	\EEX
	
	\BEX
	Time-average of fBm $Y$ with parameter $\alpha \in(0,2]$ is a Gaussian process defined as
	$$
	Y(t)=\sqrt{\alpha+2} \frac{1}{t} \int_0^t B_\alpha(s) d s .
	$$
	In view of \cite{tabis} $\beta=1$, $Y \in \mathbf{S}(\alpha,2,1)$ and $R=\alpha/2+1$, if $\alpha>1$, $R=2$, if $\alpha=1$. Let $\gamma=\beta/2=1/2$. 
	If $\alpha>1$, then \eqref{mainasymp} holds with $L_u=Lu^{-2}$, $p=0$ and $c=e^{-2^{\frac{1}{2\alpha}}dL^{1/2}-2^{\frac{1}{\alpha}}\left(\frac{\alpha+2}{4}\right)L}$. If $\alpha =1$, then \eqref{mainasymp} holds with $p=0$ and $c=\MB_{Y}^{h}
	\left(\frac{L}{\sqrt{2}},[0,\infty)\right)$, where $h(t)=t+\sqrt{2}dt^{1/2}$.
	\EEX
	
	\BEX
	Dual fBm $Y$ with parameter $\alpha \in (0,2)$ a centered
	Gaussian process with covariance function 
	$$
	R_Y(t,s)= \frac{t^{\alpha}s+s^{\alpha}t}{t+s}.
	$$
	In view of \cite{tabis} $Y \in \mathbf{S}(\alpha,2,\alpha/2)$ and $\beta=1$ and $R=2$, if $\alpha>1$, $R=3$, if $\alpha=1$. Let $\gamma=\beta/2$. 
	If $\alpha>1$, then \eqref{mainasymp} holds with $L_u=Lu^{-2}$, $p=0$, $$c=\exp\left(-2^{\frac{1}{2\alpha}-1/4}dL^{1/2}-2^{\frac{1}{\alpha}-1/2}L\right).$$
	If $\alpha =1$, then \eqref{mainasymp} holds with $L_u=Lu^{-2}$, $p=0$, $c=\MB_{Y}^{h}
	\left(\frac{L}{\sqrt{2}},[0,\infty)\right)$, where $h(t)=3t+\sqrt{2}dt^{1/2}$.
	\EEX

	\section{Auxiliary lemmas}
	We have to shorten the segment $[0,T]$, since 
	for $t\gg u^{\min\left(-1/\hat{\gamma},-2/\hat{\beta},-2/\kappa\right)+p}$ in \cite[Lm 4.1]{debicki2023sojourn} we are not able to get a limiting process.
	
	Recall that $\hat{\beta}=\frac{\beta \kappa}{\alpha},\;\hat{\gamma}=\frac{\gamma \kappa}{\alpha},\;\widehat{Y}\left(t\right)=Y\left(t^{\kappa/\alpha}\right)$. Under the notation of \Cref{mainth} denote $\sigma_X(t) = \mathrm{Var}(X(t))$. For a random variable $Z$, we set $\overline{Z}=\frac{Z}{\sqrt{\mathrm{Var}\left(Z\right)}}$
	\\if $\mathrm{Var}\left(Z\right)>0$.
	
	\BEL \label{Lem25}
	Suppose that $ \liminf_{u \to \infty} f(u)/u, \limsup_{u \to \infty} f(u)/u \in (0,\infty)$. There exist absolute constants $\delta, F, G > 0 $ such that
	
	\begin{align*} 
	\pk*{ 
		\begin{aligned}
		&
		\sup_{t \in [A,A+T]u^{-2/\kappa}} \overline{X}(t^{\kappa/\alpha}) > f(u) ,
		\\&
		\sup_{t \in [t_0,t_0+T]u^{-2/\kappa}} \overline{X}(t^{\kappa/\alpha}) > f(u) 
		\end{aligned}
	} 
	\leq \; F T^2 \exp \left( - G (t_0-(A+T))^\kappa \right) \Psi(f(u)) 
	\end{align*}
	for all $T \geq 1$, $ t_0 > A + T > A > 0 $ and any $ u \geq u_0 = (4\delta)^{-\kappa/2} (t_0 + T)^{\kappa/2} $, i.e. $ (t_0 + T)u^{-2/\kappa} \leq 4\delta $.
	\\
	\EEL
	
	\BEL \label{Lem27} 
	There exist absolute constants $\delta, F, G > 0 $ such that 
	\begin{align*}
	\pk*{
		\begin{aligned}
			&
			\sup_{t \in [A,A+T]u^{-2/\kappa}} \overline{X}(t^{\kappa/\alpha}) > u , 
			\\
			&
			\sup_{t \in [A+T,A+2T]u^{-2/\kappa}} \overline{X}(t^{\kappa/\alpha}) > u
		\end{aligned}
	 }
	\leq \; F  \left( T^2 \exp \left( - G \sqrt{T^\kappa} \right) + \sqrt{T} \right) \Psi(u) \;  
	\end{align*}
	for all $ A> 0 $, $T\geq 1$ and any $ u \geq u_0 = (4\delta)^{-\kappa/2} (A+2T)^{\kappa/2} $, i.e. $ (A+2T)u^{-2/\kappa} \leq 4\delta$.
	\\
	\EEL
	
	\BEL\label{betacut}
	\begin{enumerate}[(i)]
		\item If  $\alpha \leq \beta$, then for all small enough $\delta>0$
		\bqn{ \lim_{M \to \infty} \underset{u \to \infty}{\limsup}\; \frac{1}{u^{2/\kappa-2/\hat{\beta}}\Psi(u)}\pk*{ \underset{t \in [Mu^{-2/\hat{\beta}},\delta]}{\sup} X(t^{\kappa/\alpha})>u}  &=& 0;
		}
		
		\item If  $\alpha>\beta$, then for all small enough $\delta>0$  
		\bqn{ \lim_{M \to \infty} \underset{u \to \infty}{\limsup}\; \frac{1}{\Psi(u)} \pk*{ \underset{t \in [Mu^{-2/\hat{\beta}},\delta]}{\sup} X(t^{\kappa/\alpha})>u} &=& 0.
		}
	\end{enumerate}
	\EEL
	
	\BEL\label{gammacut}
	\begin{enumerate}[i)] 
		\item  If $\alpha \leq 2\gamma$, then for all small enough $\delta>0$
		\bqn{ \limit{M}\underset{u \to \infty}{\limsup}\; \frac{1}{u^{2/\kappa-1/\hat{\gamma}}\Psi(u)}\pk*{\underset{t \in [Mu^{-1/\hat{\gamma}},\delta]}{\sup}\left( \overline{X}(t^{\kappa/\alpha})-u-dt^{\hat{\gamma}}\right)>0} = 0;
		} 
		
		\item  If $\alpha > 2\gamma$, then for all small enough $\delta>0$
		\bqn{ \limit{M}\underset{u \to \infty}{\limsup}\; \frac{1}{\Psi(u)}\pk*{\underset{t \in [Mu^{-1/\hat{\gamma}},\delta]}{\sup}
				\left(
				\overline{X}(t^{\kappa/\alpha})-u-dt^{\hat{\gamma}}
				\right)
				>0} = 0. 
		}
	\end{enumerate}
	\EEL
	
	\BEL\label{kappacut}
	If $\alpha>\kappa$, let $\epsilon \in [0,2/\kappa-2/\alpha)$,
	then for $L_u=Lu^{-2/\alpha -\epsilon}$ with $L>0$ for all small enough $\delta>0$
	\BQNY
	\lim_{M \to \infty}
	\underset{u \to \infty}{\limsup} \frac{1}{\Psi\left(u\right)u^{\frac{\epsilon\alpha}{\alpha-\kappa}}}\pk*{ \int\limits_{Mu^{-2/\kappa+\frac{\epsilon\alpha}{\alpha-\kappa}}}^{\delta} \mathbb{I}\left(\overline{X}(t^{\kappa/\alpha})>u\right)t^{\kappa/\alpha-1} dt > L_u} = 0.
	\EQNY
	\EEL
	
	\BEL\label{smallcut}
	For each $p \in (0,2/\kappa)$ we have 
	\BQNY
	\lim_{M \to \infty}
	\underset{u \to \infty}{\limsup}\; \frac{1}{u^{p}\Psi\left(u\right)}\pk*{\underset{t \in [0,\frac{1}{M}u^{-2/\kappa+p}]}{\sup} \overline{X}(t^{\kappa/\alpha})>u} = 0.
	\EQNY
	\EEL
	
	Next, we will need a Pickands lemma for the sojourn functional (which is a slightly modified version of \cite[Lm 4.1]{debicki2023sojourn}).
	
	Consider
	$\xi_{u,j}\left(t\right),  t\in \KK_1, \ j\in S_u, \ {u\geq 0}$
	a family of centered Gaussian random fields (GRF's) with continuous sample paths and variance function $\sigma_{u,j}^2$. \\
	Suppose in the following that $S_u$ is a countable set for all $u$ large.\\ 
	For simplicity in the following we assume that $0\in E_1$.

	As in \cite{debicki2023sojourn} we shall impose the following assumptions:
	\begin{enumerate}[C1]
		\item\label{C0}   $\{g_{u,j}, j\in S_u\}$ is a sequence of deterministic functions of $u$ satisfying
		\BQNY
		\lim_{u\to\IF}\inf_{j\in S_u}g_{u,j}=\IF.
		\EQNY
		\item\label{C1} $\mathrm{Var}\left(\xi_{u,j}\left(\vk{0}\right)\right)=1$ for all large $u$ and any $j\in S_u$ and
		there exists some bounded continuous function $h$ on $ \KK_1$ such that
		\BQNY\label{assump-cova-field}
		\lim_{u\to\IF}\sup_{\vk{s}\in E_1 ,j\in S_u}\left|g_{u,j}^2\LT( 1-\sigma_{u,j}\left(\vk{s}\right) \RT) - h\left(\vk{s}\right)\right| =0.
		\EQNY	
		\item\label{C2} There exists a centered GRF  $\zeta\left(\vk{s}\right),\vk{s}\in \mathbb{R}^{k}$ with a.s. continuous sample paths such that 
		\BQN\label{C21}
		\lim_{u\to\IF}\sup_{s, s'\in E_1, j\in S_u}\abs{g_{u,j}^2\left(\mathrm{Var}\left(\overline\xi_{u,j}\left(\vk{s}\right)-\overline\xi_{u,j}\left(\vk{s}'\right)\right)\right) - 2\mathrm{Var}\left(\zeta\left(s\right)-\zeta\left(s'\right)\right)} =0.
		\EQN
		\item\label{C3} There exist positive constants $C, \nu, u_0$ such that
		\BQNY\label{assump-holder-field}
		\sup_{j\in S_u} g_{u,j}^2\mathrm{Var}\left(\overline\xi_{u,j}\left(\vk{s}\right)-\overline\xi_{u,j}\left(\vk{s}'\right)\right) \leq C \norm{\vk{s}-\vk{s}'}^\nu
		\EQNY
		holds for all $\vk{s},\vk{s}'\in \KK_1 , u\geq u_0$.
	\end{enumerate}

	Denote by $C\left(E_{i}\right), i=1,2$ the Banach space of  all continuous functions $f: E_i \mapsto \R$,  with $ E_{i}\subset\mathbb{R}^{k_i}, k_i\geq 1, i=1,2$ being  compact rectangles equipped with the sup-norm. 
	
	Let $\Gamma: C\left(E_{1}\right)\rw C\left(E_{2}\right)$ be a continuous functional satisfying
	\begin{enumerate}[F1)]
		\item \label{eF1}For any $f\in C\left(E_{1}\right)$, and $a>0, b\in\mathbb{R}$, $\Gamma\left(af+b\right)=a\Gamma\left(f\right)+b$;
		\item \label{eF2} There exists $c>0$ such that 
		$$\sup_{t\in E_{2}}\Gamma\left(f\right)\left(t\right)\leq c\sup_{s\in E_{1}}f\left(s\right), \ \ \forall f\in C\left(E_{1}\right).$$
	\end{enumerate}
	
	Define below for given $x\ge 0$ and for given
	positive $\sigma$-finite measure $\eta$ on $E_2$
	\BQNY
	&&\MB^{\Gamma, h,\eta}_{\zeta}\left( x, E_2\right)= \int_\R  \pk*{ \int_{ E_{2} } \mathbb{I}\left(\Gamma\left(\sqrt{2}\zeta-\mathrm{Var}\left(\zeta\right)-h\right)\left(t\right) + y >0\right)  \eta(dt) >x } e^{-y} dy.
	\EQNY

	\BEL\label{the-weak-conv}
	Let $\{\xi_{u,j}\left(\vk{s}\right),\vk{s}\in E_{1} ,j\in S_u, {u\geq 0}\}$ be a family of centered GRF's defined as above satisfying {\Cref{C0}-\Cref{C3}} and let
	$\Gamma$  satisfy  {\bf F1-F2}.
	Let $\eta$ be a positive $\sigma$-finite measure on
	$E_2$ being equivalent with the Lebesgues measure on $E_2$. If for  all large
	$u$ and all $j \in S_u$
	$$\pk*{ \sup_{t\in E_{2}} \Gamma\left(\xi_{u,j}\right)\left(t\right)> g_{u,j}} > 0,$$
	then for all $x\in[0,\eta\left(E_{2}\right))$
	\BQN\label{con-uni-con}
	\lim_{u\to\IF}\sup_{j\in S_u} \ABs{ \frac{\pk*{ \int_{ E_{2} } \mathbb{I}\LT(\Gamma\left(\xi_{u,j}\right)\left(\vk{t}\right)>g_{u,j}\RT) \eta\left(dt\right) >x} } {\Psi\left(g_{u,j}\right)} -\MB^{\Gamma, h,\eta}_{\zeta}\left( x, E_2\right) } = 0, 
	\EQN
	and the constant $\MB^{\Gamma, h,\eta}_{\zeta}\left( x, E_2\right)$ is continuous at  $x\in \left(0,\eta\left(E_{2}\right)\right)$.
	
	The convergence \eqref{con-uni-con} is uniform in $x$ in compact subsets of $\left(0,\eta\left(E_2\right)\right)$.
	\EEL
	If $E_1=E_2$ and $\Gamma$ is the identity functional, that is, $\Gamma(f)=f$ for each $f \in C(E_1)$, then we suppress the superscript $\Gamma$ and write simply $\MB^{h,\eta}_{\zeta}$ instead of $\MB^{\Gamma,h,\eta}_{\zeta}$.
	
	
	For a given level $u\inr$ define the excursion set of $X$ above the level $u$ by
	$$ A_u\left(X\right)= \{t\in E: X\left(t\right) >u \}.$$
	
	Consider the following convergence
	\begin{equation}\label{shti} 
		\limit{u} \pk*{ \mu\left(A_u\left(X\right)\right) >\vv\left(u\right)z\Bigl\lvert \sup_{t\in E} X\left(t\right)>u}= \EH{\bar F}\left(z\right)
	\end{equation}
	
	\def\Esu{E\left(u,n\right)}
	\BT\label{th1} 
	Let $E_u, u>0$  be compact set of $\R^k$ such that $\limit{u}\pk*{ \sup_{t\in E_u} X\left(t\right)>u}=0 $. 
	Suppose that there exist collections of Lebesgue measurable
	disjoint compact sets $ I_{k}\left(u,n\right), k\in K_{u,n}$ with $K_{u,n}$ non-empty countable index sets    such that
	$$
	E\left(u,n\right)= \bigcup_{k\in K_{u,n}}I_{k}\left(u,n\right)\subset E_u,
	$$
	then \eqref{shti} holds with $E=E_u,\;\mu=\mu_u$ if the measures $\mu_u$ have no atoms and the following three conditions are satisfied:
	
	A1) (Reduction to relevant sets)
	$$ \lim_{n\rw\IF}\limsup_{u\rw\IF}\frac{\pk*{ \sup_{t\in E_u\setminus\Esu} X\left(t\right)>u}}{\pk*{ \sup_{t\in \Esu} X\left(t\right)>u}}=0.$$
	A2) (Uniform single-sum approximation) There exists $v\left(u\right)>0$ and $\bar F_{u,n,k}, n\ge 1$ such that
	\begin{equation}\label{Pickands}
		\limit{u}\sup_{k\in K_{u,n}} \ABs{
			\frac{\pk*{ \mu_u\left(\{ t\in I_{k}\left(u,n\right): X\left(t\right) >u \}\right)> \vv\left(u\right)x}}
			{\pk*{ \sup_{ t\in I_{k}\left(u,n\right)} X\left(t\right)>u} }- \bar F_{u,n,k}\left(x\right)}= 0,
	\end{equation}
	where $x\geq 0,\ n\ge 1$,
	and for all $x\geq 0$
	\BQN\label{constant}
	\lim_{n\to\IF}\limsup_{u \to \IF} \left|\frac
	{\sum_{k\in K_{u,n}} \pk*{ \sup_{t\in I_{k}\left(u,n\right)} X\left(t\right)>u}\bar{F}_{u,n,k}\left(x\right)}
	{\sum_{k\in K_{u,n}} \pk*{ \sup_{t\in I_{k}\left(u,n\right)} X\left(t\right)>u}}-\bar{F}\left(x\right) \right|=0,
	\EQN
	where $\bar{F}\left(x\right) \in (0,1]$.
	
	A3) (Double-sum negligibility) For all large $n$ and large $u$,  $\sharp K_{u,n}\geq 2$
	and 
	$$\lim_{n\rw\IF}\limsup_{u\rw\IF} \frac
	{\sum_{ i \neq j, i,j\in K_{u,n}} \pk*{ \sup_{t\in I_{i}\left(u,n\right)} X\left(t\right)>u,
			\sup_{t\in I_j\left(u,n\right)} X\left(t\right)>u}}{\sum_{k\in K_{u,n}} \pk*{ \sup_{t\in I_{k}\left(u,n\right)} X\left(t\right)>u}}= 0.
	$$
	
	\ET

	\section{Proofs}
	Hereafter, $Q_i, i\in \N$ are some positive constants which might be different from line to line and $f\left(u,n\right)\sim g\left(u\right), u\to\infty, n\to\infty$ means that $$\lim_{n\to\infty}\lim_{u\to\infty}\frac{f\left(u,n\right)}{g\left(u\right)}=1.$$
	
	We shall give first the proofs of the auxiliary results and then present the proofs of the main results.
	
	\subsection{Proofs of auxiliary results}
	Remark that by the Borell inequality (see, e.g., \cite{Pit96})
	\BQN
	\label{deltacut}
	\lefteqn{\pk*{\int\limits_{0}^{T} \mathbb{I}\left(X\left(t\right)-dt^{\gamma}>u\right) dt >\left(\frac{\kappa}{\alpha}\right)L_u} }\notag \\
	&=&
	\pk*{\int\limits_{0}^{T^{\alpha/\kappa}} \mathbb{I}\left(\widehat{X}\left(t\right)-dt^{\hat{\gamma}}>u\right) t^{\kappa/\alpha-1}dt >L_u}\notag
	\\ &=&
	\pk*{\int\limits_{0}^{\delta} \mathbb{I}\left(\frac{\widehat{X}\left(t\right)}{1+\frac{dt^{\hat{\gamma}}}{u}}>u\right) t^{\kappa/\alpha-1}dt >L_u}+o\left(\Psi(u)\right)\notag
	\\
	&=&
	\pk*{\int\limits_{0}^{\delta} \mathbb{I}\left(\overline{X}(t^{\kappa/\alpha})>\frac{u\left(1+\frac{dt^{\hat{\gamma}}}{u}\right)}{\sigma_{\widehat{X}}\left(t\right)}\right) t^{\kappa/\alpha-1}dt >L_u}+o\left(\Psi(u)\right).
	\EQN
	as $u \to \infty$ for an arbitrary fixed $\delta$.
	
	Denote $c_1 = \underset{x \in [0,1]}{\inf}\frac{V_{\widehat{Y}}\left(1,x\right)}{|1-x|^{\kappa}},\;c_2 = \underset{x \in [0,1]}{\sup}\frac{V_{\widehat{Y}}\left(1,x\right)}{|1-x|^{\kappa}}$. \\
	In view of S2 we have 
	$$0 < c_1,c_2 < \infty.$$
	Pick $\delta$ positive such that $4\delta<\left(\frac{1}{8 a c_2}\right)^{1/\kappa}$ and for all $t,s \in [0,4\delta]$ we have $$1-r_{\widehat{X}}\left(t,s\right)\in \left[\frac{a}{2} V_{\widehat{Y}}\left(t,s\right),2a V_{\widehat{Y}}\left(t,s\right)\right].$$
	
	Then 
	\BQN
	\label{bddcorr}
	1-r_{\widehat{X}}\left(t,s\right) \in \left[\frac{ac_1}{2} |t-s|^{\kappa},2ac_2 |t-s|^{\kappa}\right] \text{ for all } t,s \in [0,4\delta].
	\EQN
	
	We claim that for such $\delta>0$ the conclusions of Lemmas \ref{Lem25}-\ref{smallcut} hold.

	\begin{proof}[Proof of \Cref{Lem25}]
		The proof is analogous to the proof of \cite[Lm 3.12]{tabis}. But for clarity we provide a complete proof below.	
		
		Let $ u_0 = (4\delta)^{-\kappa/2} (t_0 + T)^{\kappa/2} $ and $ \left\{ Z_u(t_1,t_2) : (t_1,t_2) \in [A,A+T]\times[t_0,t_0+T] \right\} $, where 
		\\$ Z_u(t_1,t_2) = \overline{X}((t_1 u^{-2/\kappa})^{\kappa/\alpha}) + \overline{X}((t_2 u^{-2/\kappa})^{\kappa/\alpha}) $. Note that
		\BQN &&\pk*{ \sup_{t \in [A,A+T]u^{-2/\kappa}} \overline{X}(t^{\kappa/\alpha}) > f(u) , \sup_{t \in [t_0,t_0+T]u^{-2/\kappa}} \overline{X}(t^{\kappa/\alpha}) > f(u) } \notag \\
		&&\label{lab34} \leq \; \pk*  { \sup_{(t_1,t_2) \in [A,A+T]\times[t_0,t_0+T]} Z_u(t_1,t_2) > 2f(u) } . 
		\EQN
		Hence, using the fact that $x \leq 1-e^{-2x}$ for $x \in [0,1/2]$, by \eqref{bddcorr}, as $$4ac_2u^{-2}|t_2-t_1|^{\kappa}\leq 4ac_2 (u^{-2/\kappa}(t_0+T))^{\kappa} \leq 4ac_2 (4\delta)^{\kappa} < 1/2,$$ we obtain
		\begin{eqnarray}
			\nonumber
			V_{\overline{X}} \! \left((t_1u^{-2/\kappa})^{\kappa/\alpha},(t_2u^{-2/\kappa})^{\kappa/\alpha}\right) & = & 2  \left( 1 - r_{\widehat{X}}(t_1 u^{-2/\kappa},t_2 u^{-2/\kappa}) \right)
			\\ 
			\label{lab33}
			&\geq& \; ac_1 u^{-2}|t_2-t_1|^{\kappa} \; ; 
			\\ \nonumber  
			V_{\overline{X}} \! \left((t_1u^{-2/\kappa})^{\kappa/\alpha},(t_2u^{-2/\kappa})^{\kappa/\alpha}\right) & = & 2  \left( 1 - r_{\widehat{X}}(t_1 u^{-2/\kappa},t_2 u^{-2/\kappa}) \right)
			\\
			\label{lab37}
			&\leq&
			\left( 1 - \exp \left(-8ac_2 u^{-2}|t_2-t_1|^{\kappa}\right) \right)
		\end{eqnarray}
		for all $ t_1,t_2 \leq t_0+T $. Since 
		\begin{eqnarray*}
			\sigma_{Z_u}^2(t_1,t_2) & = & 2  \left( 1 + r_{\widehat{X}}(t_1 u^{-2/\kappa},t_2 u^{-2/\kappa}) \right) \\[1ex]
			& = & 4 - V_{\overline{X}} \! \left((t_1u^{-2/\kappa})^{\kappa/\alpha},(t_2u^{-2/\kappa})^{\kappa/\alpha}\right) , 
		\end{eqnarray*}
		then from (\ref{lab33}) for any $ (t_1,t_2) \in [A,A+T]\times[t_0,t_0+T] $,
		\begin{equation} \label{lab35} 2 \; \leq \; \sigma_{Z_u}^2(t_1,t_2) \; \leq \; 4 - ac_1 u^{-2}(t_0-(A+T))^{\kappa} \; . 
		\end{equation}
		
		Now observe that
		\BQN 
		&&
		\pk*  { \sup_{(t,s) \in [A,A+T]\times[t_0,t_0+T]} Z_u(t_1,t_2) > 2f(u) } \notag\\
		&&
		\label{lab39} \leq \; \pk*  { \sup_{(t_1,t_2) \in [A,A+T]\times[t_0,t_0+T]} \overline{Z}_u(t_1,t_2) > \frac{2f(u)}{\sqrt{ 4 - ac_1 u^{-2}(t_0-(A+T))^{\kappa}}} } . 
		\EQN
		Remark that 
		\BQN
		\label{lab13}
		V_{\overline{Z}_u}((t_1,t_2),(s_1,s_2)) 
		=
		\frac{V_{Z_u}((t_1,t_2),(s_1,s_2))-(\sigma_{Z_u}(t_1,t_2)-\sigma_{Z_u}(s_1,s_2))^2}{\sigma_{Z_u}(t_1,t_2)\sigma_{Z_u}(s_1,s_2)}.
		\EQN
		Note that (using (\ref{lab13}) and (\ref{lab35})) for any $ (t_1,t_2), (s_1,s_2) \in [A,A+T]\times[t_0,t_0+T] $, we have
		\BQN &&V_{\overline{Z}_u}((t_1,t_2),(s_1,s_2)) \; \leq \; \frac{V_{Z_u}((t_1,t_2),(s_1,s_2))}{\sigma_{Z_u}(t_1,t_2) \sigma_{Z_u}(s_1,s_2)} \notag \\
		&&\leq \; \frac{1}{2} \: \mathbb{E}  \Bigg\{ \Big( \left( \overline{X}((t_1u^{-2/\kappa})^{\kappa/\alpha}) - \overline{X}((s_1u^{-2/\kappa})^{\kappa/\alpha}) \right)
		\notag
		\\
		&&\quad
		+ \left( \overline{X}((t_2u^{-2/\kappa})^{\kappa/\alpha}) - \overline{X}((s_2u^{-2/\kappa})^{\kappa/\alpha}) \right)\Big) ^2\Bigg\} \notag \\
		\label{lab36} &&\leq \; V_{\overline{X}} \! \left((t_1 u^{-2/\kappa})^{\kappa/\alpha},(s_1 u^{-2/\kappa})^{\kappa/\alpha}\right) + V_{\overline{X}} \! \left((t_2 u^{-2/\kappa})^{\kappa/\alpha},(s_2 u^{-2/\kappa})^{\kappa/\alpha}\right) 
		\\
		\label{lab38} &&\leq \;  \left( 1 - \exp \left(-8ac_2 u^{-2}|t_1-s_1|^{\kappa}\right) \right) + \left( 1 - \exp \left(-8ac_2 u^{-2}|t_2-s_2|^{\kappa}\right) \right) , 
		\EQN
		where (\ref{lab36}) follows from the inequality $ (x+y)^2 \leq 2(x^2+y^2) $ and (\ref{lab38}) follows from (\ref{lab37}). 
		
		Consider two independent, identically distributed centered stationary Gaussian processes $ \{ Z_{1,u}(t_1) : t_1 \geq 0 \} $, $ \{ Z_{2,u}(t_2) : t_2 \geq 0 \} $ with $ R_{Z_{1,u}}(t_1,s_1) = \exp \left(-8ac_2 u^{-2}|t_1-s_1|^{\kappa}\right) $ and let $ \underline{Z}_u(t_1,t_2) = \frac{1}{\sqrt{2}}  \left( Z_{1,u}(t_1) + Z_{2,u}(t_2) \right) $. Hence, by (\ref{lab38}), for any $ (t_1,t_2), (s_1,s_2) \in [A,A+T]\times[t_0,t_0+T] $,
		\BQNY &&V_{\overline{Z}_u}((t_1,t_2),(s_1,s_2)) 
		\\
		&&\leq \;  \left( 1 - \exp \left(-8ac_2 u^{-2}|t_1-s_1|^{\kappa}\right) \right) + \left( 1 - \exp \left(-8ac_2 u^{-2}|t_2-s_2|^{\kappa}\right) \right) 
		\\
		&&
		= \; V_{\underline{Z}_u}((t_1,t_2),(s_1,s_2))
		\EQNY
		and due to Slepian's inequality, we obtain that (with $ u^* = \frac{2f(u)}{\sqrt{ 4 - ac_1 u^{-2}(t_0-(A+T))^{\kappa}}}$)
		\BQN 
		&&
		\pk* { \sup_{(t_1,t_2) \in [A,A+T]\times[t_0,t_0+T]} \overline{Z}_u(t_1,t_2) > u^* } \; 
		\notag
		\\
		&&
		\leq  \pk*  { \sup_{(t_1,t_2) \in [A,A+T]\times[t_0,t_0+T]} \underline{Z}_u(t_1,t_2) > u^* } \notag \\
		&&
		\label{lab40} = \; \pk*  { \sup_{(t_1,t_2) \in [0,T]^2} \underline{Z}_u(t_1,t_2) > u^* },
		\EQN
		where the last equality follows from stationarity of $ \underline{Z}_u(\cdot,\cdot) $. Now we can apply \cite[Thm 2.2]{Uniform2016} to (\ref{lab40}), which combined with (\ref{lab34}) and (\ref{lab39}) gives
		\begin{align}
			&
			\pk*{
				\begin{aligned}
					& \sup_{t \in [A,A+T]u^{-2/\kappa}} \overline{X}(t^{\kappa/\alpha}) > f(u), \\
					& \sup_{t \in [t_0,t_0+T]u^{-2/\kappa}} \overline{X}(t^{\kappa/\alpha}) > f(u)
				\end{aligned}
			}
			\notag
			\\
			&
			 \leq
			\sum_{k=0}^{[T]} \sum_{l=0}^{[T]}
			\pk*{
				\begin{aligned}
					& \sup_{t \in [A+k,A+(k+1)]u^{-2/\kappa}} \overline{X}(t^{\kappa/\alpha}) > f(u), \\
					& \sup_{t \in [t_0+l,t_0+(l+1)]u^{-2/\kappa}} \overline{X}(t^{\kappa/\alpha}) > f(u)
				\end{aligned}
			}
			\notag
			\\&
			\leq \;
			4T^2 \left( \mathcal{H}_{B_\kappa} \! \left((4ac^2c_2)^{1/\kappa}\right) \right)^2 \Psi(u^*) (1+o(1)),
			\label{lab41}
		\end{align}
		uniformly in $T \geq 1$
		as $ u \to \infty $. Since 
		\begin{eqnarray*}
			(u^*)^2 & = & \frac{4f^2(u)}{4 - ac_1 u^{-2}(t_0-(A+T))^{\kappa}} \; \geq \; f^2(u) + \frac{a c_1}{4}  \left( \frac{f(u)}{u} \right)^2 (t_0-(A+T))^{\kappa} \\[1ex]
			& \geq & f^2(u) + \frac{\underline{c}^2 ac_1}{4} (t_0-(A+T))^{\kappa} \; ,
		\end{eqnarray*}
		where $ \underline{c} > 0 $ is a constant such that $ \frac{f(u)}{u} \geq \underline{c} $ for all $ u \geq u_0 $, then, since 
		\\
		$\Psi(u) \leq \frac{\exp(-u^2/2)}{u \sqrt{2\pi}}$ for any $u>0$, 
		\begin{eqnarray}
			\nonumber \Psi(u^*) & \leq & \frac{\exp  \left( - \frac{1}{2}  \left( f^2(u) + \frac{\underline{c}^2 ac_1}{4} (t_0-(A+T))^{\kappa} \right) \right)}{\sqrt{2 \pi} \sqrt{f^2(u) + \frac{\underline{c}^2 ac_1}{4} (t_0-(A+T))^{\kappa}}} \\[1ex]
			\label{lab42} & \leq & \frac{\exp  \left( - \frac{1}{2} f^2(u) \right)}{\sqrt{2 \pi} f(u)} \exp  \left( - \frac{\underline{c}^2 ac_1}{8} (t_0-(A+T))^{\kappa} \right) .
		\end{eqnarray}
		By (\ref{lab41}) and (\ref{lab42}) we obtain that
		\BQN &&\pk*  { \sup_{t \in [A,A+T]u^{-2/\kappa}} \overline{X}(t) > f(u) , \sup_{t \in [t_0,t_0+T]u^{-2/\kappa}} \overline{X}(t) > f(u) } \notag \\
		&&
		\leq \; 4 F_1 (4c^2ac_2)^{2/\kappa}  \left( \mathcal{H}_{B_\kappa}(1) \right)^2 T^2 \exp  \left( - \frac{\underline{c}^2 ac_1}{8} (t_0-(A+T))^{\kappa} \right) \Psi(f(u)) \; , 
		\EQN
		for any $ u \geq u_0 $ and some positive constant $ F_1 $, such that $ \frac{\exp \left( - \frac{1}{2} f^2(u) \right)}{\sqrt{2 \pi} f(u)} \leq F_1 \Psi(f(u)) $ for $ u > u_0 $.
		This completes the proof with $ F = 4 F_1 (4c^2 ac_2)^{2/\kappa}  \left( \mathcal{H}_{B_\kappa}(1) \right)^2 $ and $ G = \frac{\underline{c}^2 ac_1}{8} $.	
	\end{proof}
	
	\begin{proof}[Proof of \Cref{Lem27}]
		The proof is almost the same as the proof of \cite[Lm 3.13]{tabis}, but for completeness we provide it below.	
		
		Let $ u_0 = (4\delta)^{-\kappa/2} (A+2T)^{\kappa/2} $ and $ \overline{X}_u(t) = \overline{X}((tu^{-2/\kappa})^{\kappa/\alpha}) $. We have
		\begin{align*} 
		&\pk*  { \sup_{t \in [A,A+T]} \overline{X}_u(t) > u , \sup_{t \in [A+T,A+2T]} \overline{X}_u(t) > u } \\
		& = \;
		\pk*{
		\begin{aligned}
		& \sup_{t \in [A,A+T]} \overline{X}_u(t) > u , \\& \left\{ \sup_{t \in [A+T,A+T+\sqrt{T}]} \overline{X}_u(t) > u \vee \sup_{t \in [A+T+\sqrt{T},A+2T]} \overline{X}_u(t) > u \right\} 
		\end{aligned}
		}
		\\
		& \leq \; \pk*  { \sup_{t \in [A,A+T]} \overline{X}_u(t) > u , \sup_{t \in [A+T+\sqrt{T},A+2T+\sqrt{T}]} \overline{X}_u(t) > u } \\
		& \quad \; + \; \pk*  { \sup_{t \in [A+T,A+T+\sqrt{T}]} \overline{X}_u(t) > u }  \\
		& \leq \; F_1 T^2 \exp \left( - G \sqrt{T^\kappa} \right) \Psi(u) \; + \; \pk*  { \sup_{t \in [A+T,A+T+\sqrt{T}]} \overline{X}_u(t) > u } , 
	\end{align*}
		where the last inequality follows from Lemma \ref{Lem25} with $ t_0 = A + T + \sqrt{T} $. Let $\widetilde{X}$ be a stationary centered GRF with $R_{\widetilde{X}}(t,s)=\exp(-8ac_2 |t-s|^{\kappa})$. Due to Slepian's inequality and the stationarity of $\widetilde{X}$
		\BQNY
		\pk*  { \sup_{t \in [A+T,A+T+\sqrt{T}]} \overline{X}_u(t) > u }
		&\leq&
		\pk*  { \sup_{t \in [A+T,A+T+\sqrt{T}]u^{-2/\kappa}} \widetilde{X}(t) > u }
		\\
		&\leq& 
		\sum_{k=0}^{[\sqrt{T}]}
		\pk*  { \sup_{t \in [A+T+k,A+T+(k+1)]u^{-2/\kappa}} \widetilde{X}(t) > u }
		\\
		&=&
		2 \sqrt{T}
		\pk*  { \sup_{t \in [0,1]u^{-2/\kappa}} \widetilde{X}(t) > u }
		\EQNY
		We can apply \cite[Thm 2.2]{Uniform2016} to the right side of the inequality above. We obtain that for sufficiently large $ u \geq u_0 $,  
		\BQN &&\pk*  { \sup_{t \in [A,A+T]} \overline{X}_u(t) > u , \sup_{t \in [A+T,A+2T]} \overline{X}_u(t) > u } \notag\\
		&&  \leq \; F_1 T^2 \exp \left( - G \sqrt{T^\kappa} \right) \Psi(u) + 4 \mathcal{H}_{B_\kappa}(a^{1/\kappa}) \sqrt{T} \Psi(u) (1+o(1)) \notag \\
		&& \leq \; F  \left( T^2 \exp \left( - G \sqrt{T^\kappa} \right) + \sqrt{T} \right) \Psi(u) \notag 
		\EQN
		for some constant $ F > 0 $. This completes the proof.	
	\end{proof}
	
	\begin{proof}[Proof of Lemma \ref{betacut}]
		The proof could be obtained by a slight modification of the proof of \cite[Thm 4.2]{tabis}. But for clarity we provide a complete proof below. Fix arbitrary $\epsilon \in (0,b)$. Setting  $q=\min\left(-2/\hat{\beta},-2/\kappa\right)$,
		$M_u= M u^{-2/\hat{\beta}}$ and $N = [\frac{\delta-M_u}{u^q}]$, if a suitably small $\delta$ is chosen, we obtain
		
		\BQNY
		&&\pk*{\underset{t \in [M_u,\delta]}{\sup} \widehat{X}(t)
		 >u}
	 	\\
	 	&&\leq 
		\sum\limits_{k=0}^{N} \pk*{\underset{t \in [M_u+ku^{q},M_u+(k+1)u^{q}] }{\sup} \overline{X}(t^{\kappa/\alpha}) > u (1+(b-\epsilon)(M_u+ku^{q})^{\hat{\beta}})}.
		\EQNY
		Remark that by (\ref{bddcorr}), if a suitably small $\delta$ is chosen, then 
		$$1-r_{\widehat{X}}(t,s) \leq 1-\exp\left(-4ac_2|t-s|^{\kappa}\right), \quad \forall t,s \in [0,4\delta].$$
		Consequently, we can apply the Slepian inequality. More specifically, let
		\\ $\{Y(t), t \geq 0\}$ be a centered Gaussian process with covariance function $R_{Y}(t,t')=\exp\left(-4ac_2|t-t'|^{\kappa}\right)$.
		
		Then we have
		\BQNY
		&&\sum\limits_{k=0}^{N} \pk*{\underset{t \in [M_u+ku^{q},
				M_u+(k+1)u^{q}] }{\sup} \overline{X}(t^{\kappa/\alpha}) > 
			u (1+(b-\epsilon)(M_u+ku^{q})^{\hat{\beta}})} \\ &&\leq
		\sum\limits_{k=0}^{N} \pk*{\underset{t \in [M_u+ku^{q},M_u+(k+1)u^{q}] }{\sup} Y(t) > 
			u (1+(b-\epsilon)(M_u+ku^{q})^{\hat{\beta}})}\\&&=
		\sum\limits_{k=0}^{N} \pk*{\underset{t \in [0,u^{q}] }{\sup} Y(t) > u (1+(b-\epsilon)(M_u+ku^{q})^{\hat{\beta}})}.
		\EQNY
		By \cite[Thm 2.2]{Uniform2016} (which gives an asymptotic approximation uniform in $k$) for some $C>0$, where $C$ is independent of $M$, the above expression is equal to
		\BQNY
		(1+o(1))C\sum\limits_{k=0}^{\infty} \Psi(u (1+(b-\epsilon)(Mu^{-2/\hat{\beta}}+ku^{q})^{\hat{\beta}})).
		\EQNY
		Since there exists $u_0>0$ such that for all $u \geq u_0$ and $c \geq 1$ we have 
		\bqn{\label{nacht} 
			\frac{\Psi(cu)}{\Psi(u)}
			\leq
			2
			\frac{(cu \sqrt{2\pi})^{-1}e^{-c^2 u^2/2}}{(u \sqrt{2\pi})^{-1}e^{- u^2/2}}=
			\frac{2}{c} \exp\left(-u^2(c^2-1)/2 \right)\leq
			2 \exp\left(-(c-1)u^2\right),
		}	
		we obtain
		\BQNY
		&&
		\sum\limits_{k=0}^{\infty} \Psi(u (1+(b-\epsilon)(M_u+
		ku^{q})^{\hat{\beta}}))
		\\
		&& \leq 
		\Psi(u)u^{-q-2/\hat{\beta}}\sum\limits_{k=0}^{\infty} \exp\left(-(b-\epsilon)(M+ku^{q+2/\hat{\beta}})^{\hat{\beta}}\right)u^{q+2/\hat{\beta}}\\
		&&
		\leq\Psi(u)u^{-q-2/\hat{\beta}}
		\left(
		\exp
		\left(-(b-\epsilon)M^{\hat{\beta}}\right)
		+
		\int\limits_{0}^{\infty}
		\exp
		\left(-(b-\epsilon)(M+x)^{\hat{\beta}}\right) dx
		\right)
		\EQNY
		as  $u \to \infty$.
		Since further 
		$\limit{M}\int\limits_{0}^{\infty}\exp\left(-(b-\epsilon)(M+x)^{\hat{\beta}}\right) dx = 0$, then the claim follows.
	\end{proof}
	
	\begin{proof}[Proof of Lemma \ref{gammacut}]
		The proof is almost the same as the proof of Lemma \ref{betacut}, but for completeness we provide it below.
		
		Fix arbitrary $\epsilon \in (0,d)$. Setting  $q=\min\left(-2/\kappa,-1/\hat{\gamma}\right)$, $N = [\frac{\delta-N_u}{u^q}]$ and $N_u=M u^{-1/\hat{\gamma}}$ we can write further 
		
		\BQNY
		&&
		\pk*{\underset{t \in [N_u,\delta]}{\sup} \widehat{X}(t)-dt^{\hat{\gamma}} >u}
		\\
		&&\leq
		\sum\limits_{k=0}^{N} \pk*{\underset{t \in [N_u+ku^{q},N_u+(k+1)u^{q}] }{\sup} \overline{X}(t^{\kappa/\alpha}) > 
			u+(d-\epsilon)
			(N_u+ku^q)^{\hat{\gamma}}}.
		\EQNY
		Applying again the Slepian inequality we have for $Y(t), t \geq 0$ a centered Gaussian process with covariance function $R_{Y}(t,t')=\exp\left(-4ac_2|t-t'|^{\kappa}\right)$
		%
		\BQNY
		&&\sum\limits_{k=0}^{N} \pk*{\underset{t \in [N_u+ku^{q},N_u +(k+1)u^{q}] }{\sup} \overline{X}(t^{\kappa/\alpha}) > 
			u+(d-\epsilon)
			(N_u+ku^q)^{\hat{\gamma}}} \\ &&\leq
		\sum\limits_{k=0}^{N} \pk*{\underset{t \in [N_u+ku^{q},N_u+(k+1)u^{q}] }{\sup} Y(t) > 
			u+(d-\epsilon)
			(N_u+ku^q)^{\hat{\gamma}}}\\&&= 
		\sum\limits_{k=0}^{N} \pk*{\underset{t \in [0,u^{q}] }{\sup} Y(t) > u+(d-\epsilon)
			(N_u+ku^q)^{\hat{\gamma}}}.
		\EQNY
		By \cite[Thm 2.2]{Uniform2016} for some $C>0$, where $C$ is independent of $M$, the above expression is equal to
		\BQNY
		(1+o(1))C\sum\limits_{k=0}^{\infty} \Psi(u+(d-\epsilon)
		(N_u+ku^q)^{\hat{\gamma}}).
		\EQNY
		Hence utilising \eqref{nacht} yields 
		\BQNY
		&&\sum\limits_{k=0}^{\infty} \Psi(u+(d-\epsilon)
		(N_u+ku^q)^{\hat{\gamma}}) 
		\\
		&&\leq
		\Psi(u)u^{-q-1/\hat{\gamma}}\sum\limits_{k=0}^{\infty} \exp\left(-(d-\epsilon)(M+ku^{q+1/\hat{\gamma}})^{\hat{\gamma}}\right)u^{q+1/\hat{\gamma}}\\&&
		\leq\Psi(u)u^{-q-1/\hat{\gamma}}
		\left(
		\exp\left(-(d-\epsilon)M^{\hat{\gamma}}\right)+
		\int\limits_{0}^{\infty}
		\exp\left(-(d-\epsilon)(M+x)^{\hat{\gamma}}\right) dx
		\right)
		\EQNY
		as $u \to \infty$.
		Since $\int\limits_{0}^{\infty}\exp\left(-(d-\epsilon)(M+x)^{\hat{\gamma}}\right) dx \to 0$ as $M \to \infty$ the proof is complete. 
	\end{proof}
	
	\begin{proof}[Proof of Lemma \ref{kappacut}]
		Remark that
		\begin{align*}
		&\pk*{ \int\limits_{Mu^{-2/\kappa+\frac{\epsilon\alpha}{\alpha-\kappa}}}^{\delta} \mathbb{I}\left(\overline{X}(t^{\kappa/\alpha})>u\right)t^{\kappa/\alpha-1} dt > L_u}\\&
		\leq\sum\limits_{j=0}^{N_u}
		\pk*{ \int\limits_{2^{j}Mu^{-2/\kappa+\frac{\epsilon\alpha}{\alpha-\kappa}}}^{2^{j+1}Mu^{-2/\kappa+\frac{\epsilon\alpha}{\alpha-\kappa}}} \mathbb{I}\left(\overline{X}(t^{\kappa/\alpha})>u\right)t^{\kappa/\alpha-1} dt > \frac{1}{2(j+1)^2}L_u}\\&
		\leq 
		\sum\limits_{j=0}^{N_u}
		\pk*{ \int\limits_{2^{j}Mu^{-2/\kappa+\frac{\epsilon\alpha}{\alpha-\kappa}}}^{2^{j+1}Mu^{-2/\kappa+\frac{\epsilon\alpha}{\alpha-\kappa}}} \mathbb{I}\left(\overline{X}(t^{\kappa/\alpha})>u\right) dt > \left(2^{j}Mu^{-2/\kappa+\frac{\epsilon\alpha}{\alpha-\kappa}}\right)^{1-\kappa/\alpha} \frac{1}{2(j+1)^2}L_u}
		\\&=
		\sum\limits_{j=0}^{N_u}
		\pk*{ \int\limits_{2^{j}Mu^{-2/\kappa+\frac{\epsilon\alpha}{\alpha-\kappa}}}^{2^{j+1}Mu^{-2/\kappa+\frac{\epsilon\alpha}{\alpha-\kappa}}} \mathbb{I}\left(\overline{X}(t^{\kappa/\alpha})>u\right) dt > \left(2^{j}M\right)^{1-\kappa/\alpha} \frac{1}{2(j+1)^2}Lu^{-2/\kappa}},
		\end{align*}
		where $N_u = \left[\log_2(\delta)-\log_2{\left(Mu^{-2/\kappa+\frac{\epsilon\alpha}{\alpha-\kappa}}\right)}\right]$.
		
		But if 
		\BQNY
		\int\limits_{2^{j}Mu^{-2/\kappa+\frac{\epsilon\alpha}{\alpha-\kappa}}}^{2^{j+1}Mu^{-2/\kappa+\frac{\epsilon\alpha}{\alpha-\kappa}}} \mathbb{I}\left(\overline{X}(t^{\kappa/\alpha})>u\right) dt > \left(2^{j}M\right)^{1-\kappa/\alpha} \frac{1}{2(j+1)^2}Lu^{-2/\kappa},
		\EQNY
		then setting 
		\BQNY
		A_k=\left\{\underset{t \in [\left(k+2^jMu^{\frac{\epsilon\alpha}{\alpha-\kappa}}\right)u^{-2/\kappa},\left(k+1+2^j Mu^{\frac{\epsilon\alpha}{\alpha-\kappa}}\right)u^{-2/\kappa}]}{\sup} \overline{X}(t^{\kappa/\alpha})>u\right\}
		\EQNY
		one obtains
		\BQNY
		\left|\left\{k \in \mathbb{Z}:\;0\leq k \leq 2^j Mu^{\frac{\epsilon\alpha}{\alpha-\kappa}},
		\,A_k \text{ holds}
		\right\}\right|>\left(2^{j}M\right)^{1-\kappa/\alpha} \frac{1}{2(j+1)^2}L,
		\EQNY
		which implies that there exist integers $0 \leq k\leq l \leq 2^j Mu^{\frac{\epsilon\alpha}{\alpha-\kappa}}$
		such that
		\\ $l-k>\left(2^{j}M\right)^{1-\kappa/\alpha} \frac{1}{2(j+1)^2}L-1$. For $M>(8/L)^{\frac{1}{1-\kappa/\alpha}}\sup_{j \geq 0} 2^{-j} (j+1)^{\frac{2}{1-\kappa/\alpha}}$ we obtain $l-k-1 > \left(2^{j}M\right)^{1-\kappa/\alpha} \frac{1}{4(j+1)^2}L$.
		
		Therefore, 
		\BQNY
		&&\pk*{ \int\limits_{2^{j}Mu^{-2/\kappa+\frac{\epsilon\alpha}{\alpha-\kappa}}}^{2^{j+1}Mu^{-2/\kappa+\frac{\epsilon\alpha}{\alpha-\kappa}}} \mathbb{I}\left(\overline{X}(t^{\kappa/\alpha})>u\right) dt > \left(2^{j}M\right)^{1-\kappa/\alpha} \frac{1}{2(j+1)^2}Lu^{-2/\kappa}}
		\\&& \leq
		\sum\limits_{0 \leq k < k +1+\left(2^{j}M\right)^{1-\kappa/\alpha} \frac{1}{4(j+1)^2}L<l \leq 2^j Mu^{\frac{\epsilon\alpha}{\alpha-\kappa}}}
		\mathbb{P}\left(A_k \cup A_l\right)
		\\&&  \leq 
		\sum\limits_{0 \leq k \leq 2^j Mu^{\frac{\epsilon\alpha}{\alpha-\kappa}}-1} \sum \limits_{l=0}^{\infty} F \exp\left(-G\left(l+\left(2^{j}M\right)^{1-\kappa/\alpha} \frac{1}{4(j+1)^2}L\right)^\kappa\right) \Psi\left(u\right)
		\\&& 
		\leq 2^j M u^{\frac{\epsilon\alpha}{\alpha-\kappa}} \sum \limits_{l=0}^{\infty} F \exp\left(-G\left(l+\left(2^{j}M\right)^{1-\kappa/\alpha} \frac{1}{4(j+1)^2}L\right)^\kappa\right) \Psi\left(u\right)
		\EQNY
		for some constants $F,G>0$, where the penultimate inequality follows from \Cref{Lem25}, if $\delta$ is chosen small enough so that the conclusion of \Cref{Lem25} holds, since one can check that
		\\ $\left(l+1+2^j Mu^{\frac{\epsilon\alpha}{\alpha-\kappa}}\right)u^{-2/\kappa} \leq 4\delta \leq \left(2ac_2\right)^{1/\kappa}$.
		
		Finally, we have
		\BQN
		\label{sojourncutbd}
		&&
		\pk*{ \int\limits_{Mu^{-2/\kappa+\frac{\epsilon\alpha}{\alpha-\kappa}}}^{\delta} \mathbb{I}\left(\overline{X}(t^{\kappa/\alpha})>u\right)t^{\kappa/\alpha-1} dt > L_u}
		\notag
		\\&& 
		\leq
		u^{\frac{\epsilon\alpha}{\alpha-\kappa}} \sum\limits_{j=0}^{\infty}\sum \limits_{l=0}^{\infty} 2^j M F \exp\left(-G\left(l+\left(2^{j}M\right)^{1-\kappa/\alpha} \frac{1}{4(j+1)^2}L\right)^\kappa\right) \Psi\left(u\right)
		\EQN
		Since the series
		\BQNY
		F\left(M\right):= \sum\limits_{j=0}^{\infty}\sum \limits_{l=0}^{\infty} 2^j M F \exp\left(-G\left(l+\left(2^{j}M\right)^{1-\kappa/\alpha} \frac{1}{4(j+1)^2}L\right)^\kappa\right)
		\EQNY
		converges for each $M$, we obtain 
		\BQN
		\label{fto0}
		F\left(2^{j'} M\right)= \sum\limits_{j=j'}^{\infty}\sum \limits_{l=0}^{\infty} 2^j M F \exp\left(-G\left(l+\left(2^{j}M\right)^{1-\kappa/\alpha} \frac{1}{4(j+1)^2}L\right)^\kappa\right) \to 0
		\EQN
		as $j' \to \infty$, but due to \eqref{sojourncutbd}
		\begin{eqnarray*}
		\underset{u \to \infty}{\limsup} \frac{1}{\Psi\left(u\right)u^{\frac{\epsilon\alpha}{\alpha-\kappa}}}\pk*{ \int\limits_{\left(2^{j'}M\right)u^{-2/\kappa+\frac{\epsilon\alpha}{\alpha-\kappa}}}^{\delta} \mathbb{I}\left(\overline{X}(t^{\kappa/\alpha})>u\right)t^{\kappa/\alpha-1} dt > L_u} 
		\leq F\left(2^{j'}M\right),
		\end{eqnarray*}
		which, combined with \eqref{fto0}, finishes the proof of Lemma \ref{kappacut}.
	\end{proof}
	
	\begin{proof}[Proof of Lemma \ref{smallcut}]
		From \cite[Thm 4.1]{tabis} it follows that there exists a constant $C>0$ such that for all $M>0$
		\BQNY
		\pk*{\underset{t \in [0,\frac{1}{M}u^{-2/\kappa+p}]}{\sup} \overline{X}(t^{\kappa/\alpha})>u}\sim\frac{C}{M}u^p \Psi\left(u\right), \quad u\to \IF,
		\EQNY
		which implies the conclusion of \Cref{smallcut}.
	\end{proof}

	\begin{proof}[Proof of Lemma \ref{the-weak-conv}]
		This is the same as \cite[Lm 4.1]{debicki2023sojourn}; the uniform convergence in $x$ on compact sets follows from the continuity of $\MB^{\Gamma, h,\eta}_{\zeta}\left( x, E_2\right)$ combined with the monotonicity of $\MB^{\Gamma, h,\eta}_{\zeta}\left( x, E_2\right)$ and\\ $\pk*{ \int_{ E_{2} } \mathbb{I}\LT(\Gamma\left(\xi_{u,j}\right)\left(\vk{t}\right)>g_{u,j}\RT) \eta\left(dt\right) >x}$ in $x$.
	\end{proof}
	\begin{proof}[Proof of Theorem \ref{th1}]
		This is a minor extension of \cite[Thm 1.1]{debicki2023sojourn} and therefore we omit the details.

	\end{proof}
	
	\subsection{Proof of main results}
	
	\begin{proof}[Proof of \Cref{last}]
		Let $X\left(t\right),\;t \in [0,T]$ be a Gaussian process with unit variance such that $r_X\left(t,s\right)=\exp\left(-V_{\widehat{Y}}\left(t,s\right)\right)$ and set  $$X_u\left(t\right)=X\left(tu^{-2/\kappa}\right).$$
		For $T>0$ let $u\left(T\right)$ be so large (we can find such an $u\left(T\right)$ by Lemma \ref{the-weak-conv}) that $u\left(T\right)^{-2/\kappa}T \to 0$ 
		\BQN 
		\label{defut}
		\left|\frac{1}{\Psi\left(u\left(T\right)\right)}
		\pk*{\int_{0}^{T} \mathbb{I}\left(X_{u\left(T\right)}\left(t\right)>u\left(T\right)\right)>x}
		-\MB_{\widehat{Y}}\left(x,[0,T]\right)\right| \leq 1.
		\EQN
		
		Next fix $S$ and $n$ two positive integers and denote
		\begin{align*}
		\Sigma\Sigma_{n,S}=
		\sum\limits_{k=n}^{n^2-1}
		\sum\limits_{l=k+1}^{n^2-1}
		\pk*{
		\begin{aligned}
		&
		\underset{[kS,\left(k+1\right)S]}{\sup} \;X_{u\left(n^2 S\right)}\left(t\right)>u\left(n^2 S\right),
		\\&
			\underset{[lS,\left(l+1\right)S]}{\sup} \;X_{u\left(n^2 S\right)}\left(t\right)>u\left(n^2 S\right)
		\end{aligned}
		}.
		\end{align*}
		We have the following bound
		\BQN
		\label{ubconst}
		&&\pk*{\int_{0}^{n^2 S} \mathbb{I}\left(X_{u\left(n^2 S\right)}\left(t\right)>u\left(n^2 S\right)\right)>x} \notag
		\\&& \leq 
		\pk*{\underset{[0,nS]}{\sup} \;X_{u\left(n^2 S\right)}\left(t\right)>u\left(n^2 S\right)}+
		\pk*{\int_{nS}^{n^2 S} \mathbb{I}\left(X_{u\left(n^2 S\right)}\left(t\right)>u\left(n^2 S\right)\right)>x} \notag
		\\&&
		\leq \pk*{\underset{[0,nS]}{\sup} \;X_{u\left(n^2 S\right)}\left(t\right)>u\left(n^2 S\right)}
		\notag
		\\
		&& 
		\quad
		+
		\sum\limits_{k=n}^{n^2-1}\pk*{\int_{kS}^{\left(k+1\right) S} \mathbb{I}\left(X_{u\left(n^2 S\right)}\left(t\right)>u\left(n^2 S\right)\right)>x}+\Sigma\Sigma_{n,S}.
		\EQN
		
		On the other hand
		\BQN
		\label{lbconst}
		\lefteqn{\pk*{\int_{0}^{n^2 S} \mathbb{I}\left(X_{u\left(n^2 S\right)}\left(t\right)>u\left(n^2 S\right)\right)>x}}\notag \\
		&\geq& \pk*{\int_{nS}^{n^2 S} \mathbb{I}\left(X_{u\left(n^2 S\right)}\left(t\right)>u\left(n^2 S\right)\right)>x}
		\notag
		\\&
		\geq& \sum\limits_{k=n}^{n^2-1}\pk*{\int_{kS}^{\left(k+1\right) S} \mathbb{I}\left(X_{u\left(n^2 S\right)}\left(t\right)>u\left(n^2 S\right)\right)>x}-\Sigma\Sigma_{n,S}.
		\EQN
		Write next for $S>0$
		$$f\left(S\right)= \underset{n \to \infty}{\limsup} \frac{\Sigma\Sigma_{n,S}}{n^2 S \Psi\left(u\left(n^2 S\right)\right)}.  $$
		Since $\lim_{n \to \infty}u(n^2 S)^{-2/\kappa} n^2 S =0$, by \Cref{Lem25} and
		\Cref{Lem27} 
		$$\limit{S}f\left(S\right) = 0.$$
		
		Consequently, letting $n \to \infty$ in (\ref{ubconst}) and (\ref{lbconst}), divided by $n^2 S$, using (\ref{defut}) for the left-hand side and the fact that (due to \cite[Thm 4.1]{tabis})
		$$\pk*{\underset{[0,nS]}{\sup} \;X_{u\left(n^2 S\right)}\left(t\right)>u\left(n^2 S\right)}=O\left(n\right)$$
		combined with Lemma \ref{the-weak-conv} for the right-hand sides we obtain
		\BQNY
		\frac{
			\MB_{B_{\kappa}\left(c^{1/\kappa}_{\widehat{Y}}\cdot\right)}\left(x,[0,S]\right)
		}
		{S}-f\left(S\right)
		&\leq& \underset{n}{\liminf}
		\frac{
			\MB_{\widehat{Y}}\left(x,[0, n^2 S]\right)
		}
		{n^2 S}
		\\
		&\leq&
		\underset{n}{\limsup}\frac{
			\MB_{\widehat{Y}}\left(x,[0, n^2 S]\right)
		}{n^2 S}
		\leq 
		\frac{\MB_{
				B_{\kappa}
				\left(c^{1/\kappa}_{\widehat{Y}} \cdot\right)}
		}
		{S}+f\left(S\right),
		\EQNY
		which implies
		\BQNY
		\frac{
			\MB_{B_{\kappa}\left(
				c^{1/\kappa}_{\widehat{Y}} \cdot\right)}\left(x,[0, S]\right)
		}
		{S}-f\left(S\right)
		&\leq& \underset{T}{\liminf}\frac{\MB_{\widehat{Y}}\left(x,[0, T]\right)}{T}
		\\
		&\leq&
		\underset{T}{\limsup}\frac{\MB_{\widehat{Y}}\left(x,[0, T]\right)}{T}
		\leq 
		\frac{
			\MB_{B_{\kappa}\left(
				c^{1/\kappa}_{\widehat{Y}} \cdot\right)}\left(x,[0, S]\right)
		}
		{S}+f\left(S\right).
		\EQNY
		Letting $S \to \infty$, we obtain 
		$\MB_{\widehat{Y}}\left(x\right)=\MB_{B_{\kappa}\left(c_{\widehat{Y}}^{1/\kappa} \cdot\right)}\left(x\right)=c_{\widehat{Y}}^{1/\kappa}\MB_{B_{\kappa}}\left(c_{\widehat{Y}}^{1/\kappa} x\right)$.
	\end{proof}

	\begin{proof}[Proof of Theorem \ref{mainth}]
		In the cases when $p\left(\alpha,\beta,\gamma,\kappa\right)>0$
		(that is, cases (i) and (ii)), we can fix $M'>0$ and apply Theorem \ref{th1} to the process $\frac{\widehat{X}\left(t\right)}{1+\frac{dt^{\hat{\gamma}}}{u}}$ with 
		$$E_u = \left[\frac{1}{M'}u^{-2/\kappa+p},M' u^{-2/\kappa+p}\right], \quad I_k\left(u,n\right)=\left[t_{u,n,k},\,t_{u,n,k+1}\right],$$
		$$t_{u,n,k}=\frac{1}{M'}u^{-2/\kappa+p}+kn\left(ac_{\widehat{Y}}\right)^{-1/\kappa}u^{-2/\kappa},$$
		where $0 \leq k \leq N\left(u,n\right)$ with
		$$N\left(u,n\right)=\left[\frac{\left(M'-\frac{1}{M'}\right)u^{p}
			\left(ac_{\widehat{Y}}\right)^{1/\kappa}}{n}\right]-1.$$
		Also, substitute $\mu_u\left(dt\right)=t^{\kappa/\alpha-1}dt, x=L, v\left(u\right)=u^{-\frac{2+p\left(\alpha-\kappa\right)}{\alpha}}$.
		Remark that in the cases when $p\left(\alpha,\beta,\gamma,\kappa\right)>0$ we have $L_u = v\left(u\right)x$.
		
		\underline{\emph{Condition A1}} 
		Remark that
		\BQNY \pk*{\underset{t \in [t_{u,n,N\left(u,n\right)},M'u^{-2/\kappa+p}]}{\sup}
			\frac{\widehat{X}\left(t\right)}{1+\frac{dt^{\hat{\gamma}}}{u}}
			>u}
		&\leq& \pk*{\underset{t \in [t_{u,n,N\left(u,n\right)},M'u^{-2/\kappa+p}]}{\sup}
			\overline{X}(t^{\kappa/\alpha})>u} 
		\\ &\leq& \pk*{\underset{t \in [t_{u,n,N\left(u,n\right)},M'u^{-2/\kappa+p}]}{\sup}
			\overline{X_2}\left(t\right)>u}
		\EQNY
		for large $u$ by \eqref{bddcorr} and
		the Slepian lemma, where $X_2$ is a centered stationary Gaussian process with unit variance such that $R_{\overline{X_2}}\left(t,s\right)=\exp\left(-4ac_2|t-s|^{\kappa}\right)$. Therefore, since $$M'u^{-2/\kappa+p}-t_{u,n,N\left(u,n\right)} \leq n\left(ac_{\widehat{Y}}\right)^{-1/\kappa}u^{-2/\kappa}$$
		by \cite[Thm 2.2]{Uniform2016} we obtain
		$$\pk*{\underset{t \in [t_{u,n,N\left(u,n\right)},M'u^{-2/\kappa+p}]}{\sup}
			\frac{\widehat{X}\left(t\right)}{1+\frac{dt^{\hat{\gamma}}}{u}}
			>u}=O\left(\Psi\left(u\right)\right),\quad u \to \infty.$$
		As we will see later
		\BQN
		\label{A1finish}
		\pk*{\sup_{t \in E\left(u,n\right)} \frac{\widehat{X}\left(t\right)}{1+\frac{dt^{\hat{\gamma}}}{u}}>u} \geq c \Psi\left(u\right)u^p 
		\EQN
		for all $u>u_n, n>n_0$, if we properly choose $u_n>0,c>0,n_0>0$. Therefore, condition A1 is proven.
		
		\underline{\emph{Condition A3}} 
		We have 
		\BQNY
		\lefteqn{
			\sum_{ i < j, i,j\in K_{u,n}} \pk*{ \sup_{t\in I_{i}(u,n)} 
				\frac{\widehat{X}(t)}{1+\frac{dt^{\hat{\gamma}}}{u}}>u,
				\sup_{t\in I_j(u,n)} \frac{\widehat{X}(t)}{1+\frac{dt^{\hat{\gamma}}}{u}}>u} }\\
		&\leq & 
		\sum_{ i < j, i,j\in K_{u,n}} \pk*{ \sup_{t\in I_{i}(u,n)} 
			\overline{X}(t^{\kappa/\alpha})>u,
			\sup_{t\in I_j(u,n)} \overline{X}(t^{\kappa/\alpha})>u}
		\\&\leq &\Sigma_{1,n}+\Sigma_{2,n},
		\EQNY
		where
		\BQNY
		\Sigma_{1,n} &=& \sum\limits_{k=0}^{N(u,n)-2}
		\sum\limits_{l=2}^{N(u,n)-k} \pk*{\sup_{t\in I_{k}(u,n)} 
			\overline{X}(t^{\kappa/\alpha})>u,
			\sup_{t\in I_{k+l}(u,n)} \overline{X}(t^{\kappa/\alpha})>u},
		\\
		\Sigma_{2,n}&=&\sum\limits_{k=0}^{N(u,n)-1}
		\pk*{\sup_{t\in I_{k}(u,n)} 
			\overline{X}(t^{\kappa/\alpha})>u,
			\sup_{t\in I_{k+1}(u,n)} \overline{X}(t^{\kappa/\alpha})>u}.
		\EQNY
		By applying \cite[Lm 4.3]{tabis} after \Cref{Lem25} for large enough $u$ and $n$ and some positive constants $F_1, G_1$ we have
		\BQNY
		\Sigma_{1,n} &\leq&
		(a c_{\hat{Y}})^{-1/\kappa} \sum\limits_{k=0}^{N(u,n)-2}
		\sum\limits_{l=2}^{N(u,n)-k} 
		F_1 n^2
		\exp\left(-G_1\left((l-1)n\right)^{\kappa}\right) \Psi(u)
		\\ &\leq&
		(a c_{\hat{Y}})^{-1/\kappa}
		\sum\limits_{k=0}^{N(u,n)-2}
		\sum\limits_{l=1}^{\infty} 
		F_1 n^2
		\exp\left(-G_1 \left(l n\right)^{\kappa}\right)\Psi(u)
		\\
		&\leq& 
		2 
		(a c_{\hat{Y}})^{-1/\kappa}
		N(u,n)  F_1 n^2
		\exp\left(-G_1 n^{\kappa}\right) \Psi(u)
		\\ &\leq& \frac{2}{n}
		M'u^p  F_1 n^2
		\exp\left(-G_1 n^{\kappa}\right) \Psi(u).
		\EQNY
		Set $t_{u,n,k,l}:=
		\frac{l}{n} t_{u,n,k+1}+
		\frac{n-l}{n} t_{u,n,k}$ and remark that $t_{u,n,k,0}=t_{u,n,k}$, $t_{u,n,k,n}=t_{u,n,k+1}$ and $t_{u,n,k,l+1}-t_{u,n,k,l}=(ac_{\widehat{Y}})^{-1/\kappa} u^{-2/\kappa}$ for all $0 \leq k \leq N(u,n)-1$ and $0 \leq l \leq n-1$. Therefore,
		\BQNY
		\limsup_{u \to \infty}
		\sup_{\substack{0 \leq k \leq N(u,n)-1,\\0 \leq l \leq n-1}}\;
		\sup_{t,s \in [t_{u,n,k,l},t_{u,n,k,{l+1}}]} 
		\ABs{\frac{1-r_{\widehat{X}}(t,s)}{|t-s|^{\kappa}}-c_{\hat{Y}}}=0.
		\EQNY
		By applying \cite[Lm 4.3]{tabis} after \Cref{Lem27} for $u$ and $n$ large enough and some positive constants $F_2,G_2$ by \Cref{the-weak-conv} with $x=0$ we have
		\BQNY
		\Sigma_{2,n} &\leq& 
		(a c_{\hat{Y}})^{-1/\kappa} N(u,n) 
		F_2 \left(n^2 \exp\left(-G_2 \sqrt{n^{\kappa}}\right)+\sqrt{n}\right)
		\Psi(u)
		\\ 
		&\leq& M'u^p \frac{1}{n}
		F_2\left(n^2
		\exp\left(-G_2\sqrt{n^{\kappa}}\right)
		+
		\sqrt{n}\right)\Psi(u)
		\EQNY
		implying further 
		\BQNY
		\underset{n \to \infty}{\lim}
		\underset{u \to  \infty}{\limsup}\frac{1}{\Psi(u)u^p} \sum_{ i \neq j,\; i,j\in K_{u,n}}\pk*{ 	\sup_{t\in I_{i}(u,n)} 
			\overline{X}(t^{\kappa/\alpha})>u,
			\sup_{t\in I_j(u,n)} \overline{X}(t^{\kappa/\alpha})>u}=0.
		\EQNY
		As we will see later from \eqref{supsinglesum}, we also have
		\BQN
		\label{singlesumbig}
		\underset{n \to \infty}{\liminf}\;
		\underset{u \to  \infty}{\liminf}\frac{1}{\Psi(u)u^p} \sum_{ i\in K_{u,n}}\pk*{\sup_{t\in I_{i}(u,n)} 
			\frac{\widehat{X}(t)}{1+\frac{dt^{\hat{\gamma}}}{u}}>u}>0
		\EQN
		and thus condition A3 follows.
		
		\underline{\emph{Condition A2}}
		Let $t = \frac{1}{M'}u^{p-2/\kappa}+\left(\tau+kn\right)
		\left(ac_{\widehat{Y}}\right)^{-1/\kappa}u^{-2/\kappa}$, $\tau \in [0,n]$.
		Then 
		\begin{align*}
		&\pk*{ \mu_u\left(\{ t\in I_{k}\left(u,n\right): \frac{\widehat{X}\left(t\right)}{1+
				\frac{dt^{\hat{\gamma}}}{u}} >u \}\right)> \vv\left(u\right)x}
		\\&=	
		\pk*{
			\int_{I\left(n,k\right)} t^{\kappa/\alpha-1}
			\mathbb{I}\left(\frac{\widehat{X}\left(t\right)}{1+
				\frac{dt^{\hat{\gamma}}}{u}}>u\right)dt>
			Lu^{-\frac{2+p\left(\alpha-\kappa\right)}{\alpha}}	
		}\\&
		=\pk*{
			\int_{[0,n]} \left(\frac{1}{M'}+
			\frac{\left(\tau+kn\right)\left(ac_{\widehat{Y}}\right)^{-1/\kappa}}{u^p}\right)^{\kappa/\alpha-1} \mathbb{I}\left(\frac{\widehat{X}\left(t\right)}{1+
				\frac{dt^{\hat{\gamma}}}{u}}>u\right)d\tau
			>L\left(ac_{\widehat{Y}}\right)^{1/\kappa}
		}\\&
		= \pk*{
			\begin{aligned}
			&
			\int_{[0,n]}  \mathbb{I}\left(\overline{X}(t^{\kappa/\alpha})> 
			\frac{
				u+
				dt_{*}^{\hat{\gamma}}
			}
			{
				\sigma_{\widehat{X}}\left(t_{**}\right)
			}
			\right)
			d\tau
			\\
			&
			>L\left(ac_{\widehat{Y}}\right)^{1/\kappa} \left(\frac{1}{M'}+
			\frac{\left(\tau^*+kn\right)\left(ac_{\widehat{Y}}\right)^{-1/\kappa}}{u^p}\right)^{1-\kappa/\alpha}
			\end{aligned}
		}
		\end{align*}
		for some $\tau^* \in [0,n]; t_{*},\, t_{**} \in I_k\left(u,n\right)$. 
		
		In turn, if $n$ is large enough, then by Lemma \ref{the-weak-conv}, since in Lemma \ref{the-weak-conv} the convergence is uniform in $x$, uniformly in $k$ we obtain (\Cref{C0}, \Cref{C1} hold trivially, \Cref{C2} follows from the definition of $c_{\widehat{Y}}$, \Cref{C3} follows from \eqref{bddcorr})
		\begin{align*}
		&
		\pk*{
			\begin{aligned}
				&
				\int_{[0,n]}  \mathbb{I}\left(\overline{X}(t^{\kappa/\alpha})> 
				\frac{
					u+
					dt_{*}^{\hat{\gamma}}
				}
				{
					\sigma_{\widehat{X}}\left(t_{**}\right)
				}
				\right)
				d\tau
				\\
				&
				>L\left(ac_{\widehat{Y}}\right)^{1/\kappa} \left(\frac{1}{M'}+
				\frac{\left(\tau^*+kn\right)\left(ac_{\widehat{Y}}\right)^{-1/\kappa}}{u^p}\right)^{1-\kappa/\alpha}
			\end{aligned}
		}
		\\&\sim\Psi\left(\frac{
			u+
			dt_{*}^{\hat{\gamma}}
		}
		{
			\sigma_{\widehat{X}}\left(t_{**}\right)
		}\right)
		\MB_{B_{\kappa}}\left(L\left(ac_{\widehat{Y}}\right)^{1/\kappa} \left(\frac{1}{M'}+
		\frac{\left(\tau^{*}+kn\right)\left(ac_{\widehat{Y}}\right)^{-1/\kappa}}{u^p}\right)^{1-\kappa/\alpha},[0,n]\right)
		\\&\sim\Psi\left(u\right)
		\MB_{B_{\kappa}}\left(L\left(ac_{\widehat{Y}}\right)^{1/\kappa} \tau_{u,n,k}^{1-\kappa/\alpha},[0,n]\right)
		\\&
		\quad 
		\times
		\exp\left(-du^{1+\left(p-2/\kappa\right)\hat{\gamma}} \tau_{u,n,k}^{\hat{\gamma}}
		-bu^{2+\left(p-2/\kappa\right)\hat{\beta}} \tau_{u,n,k}^{\hat{\beta}}\right)
		\end{align*}
		uniformly in $k$, where we define $\tau_{u,n,k}$ as $t_{u,n,k}u^{2/\kappa-p}$.
		
		By letting $x=L=0$ in the above asymptotics we obtain 
		\BQNY
		\pk*{ \sup_{ t\in I_{k}\left(u,n\right)} \frac{\widehat{X}\left(t\right)}{1+
				\frac{dt^{\hat{\gamma}}}{u}}>u}
		&\sim&
		\Psi\left(u\right)
		\MB_{B_{\kappa}}\left(0,[0,n]\right)
		\\
		&&
		\times
		\exp\left(-du^{1+\left(p-2/\kappa\right)\hat{\gamma}} \tau_{u,n,k}^{\hat{\gamma}}
		-bu^{2+\left(p-2/\kappa\right)\hat{\beta}} \tau_{u,n,k}^{\hat{\beta}}\right)
		\EQNY
		uniformly in $k$. Therefore
		$$\bar{F}_{u,n,k}\left(x\right)=\frac{
			\MB_{B_{\kappa}}\left(x\left(ac_{\widehat{Y}}\right)^{1/\kappa} \tau_{u,n,k}^{1-\kappa/\alpha},[0,n]\right)
		}{\MB_{B_{\kappa}}\left(0,[0,n]\right)}$$ 
		satisfies \eqref{Pickands}. Also we have as $u\to \IF$
		\BQN
		\label{supsinglesum}
		&&\sum_{k \in  K_{u,n}} \pk*{ \sup_{ t\in I_{k}\left(u,n\right)} \frac{\widehat{X}\left(t\right)}{1+
				\frac{dt^{\hat{\gamma}}}{u}}>u} \notag
		\\&&\sim\sum_{k \in  K_{u,n}} \Psi\left(u\right)
		\MB_{B_{\kappa}}\left(0,[0,n]\right)
		\exp\left(-du^{1+\left(p-2/\kappa\right)\hat{\gamma}} \tau_{u,n,k}^{\hat{\gamma}}
		-bu^{2+\left(p-2/\kappa\right)\hat{\beta}} \tau_{u,n,k}^{\hat{\beta}}\right) \notag
		\\&&\sim
		\frac{u^p}{n\left(ac_{\widehat{Y}}\right)^{-1/\kappa}}
		\Psi\left(u\right) \MB_{B_{\kappa}}\left(0,[0,n]\right)
		\notag
		\\
		&&
		\quad
		\times
		\sum_{k \in  K_{u,n}} \exp\left(-du^{1+\left(p-2/\kappa\right)\hat{\gamma}} \tau_{u,n,k}^{\hat{\gamma}}
		-bu^{2+\left(p-2/\kappa\right)\hat{\beta}} \tau_{u,n,k}^{\hat{\beta}}\right)\left(\tau_{u,n,k+1}-\tau_{u,n,k}\right)\notag
		\\&&
		\sim
		\frac{u^p}{n\left(ac_{\widehat{Y}}\right)^{-1/\kappa}}
		\Psi\left(u\right) \MB_{B_{\kappa}}\left(0,[0,n]\right) 
		\notag
		\\
		&&
		\quad 
		\times
		\int\limits_{\frac{1}{M'}}^{M'}
		\exp\left(-du^{1+\left(p-2/\kappa\right)\hat{\gamma}} z^{\hat{\gamma}}
		-bu^{2+\left(p-2/\kappa\right)\hat{\beta}} z^{\hat{\beta}}\right) dz
		\EQN
		and 
		\BQN
		\label{sojsinglesum}
		&&\sum_{k \in  K_{u,n}} \pk*{ \sup_{ t\in I_{k}\left(u,n\right)} \frac{\widehat{X}\left(t\right)}{1+
				\frac{dt^{\hat{\gamma}}}{u}}>u}\bar{F}_{u,n,k}\left(x\right)
		\notag
		\\&&
		\sim\sum_{k \in  K_{u,n}} \Psi\left(u\right)
		\MB_{B_{\kappa}}\left(L\left(ac_{\widehat{Y}}\right)^{1/\kappa} \tau_{u,n,k}^{1-\kappa/\alpha},[0,n]\right)
		\notag
		\\&&\quad
		\times
		\exp\left(-du^{1+\left(p-2/\kappa\right)\hat{\gamma}} \tau_{u,n,k}^{\hat{\gamma}}
		-bu^{2+\left(p-2/\kappa\right)\hat{\beta}} \tau_{u,n,k}^{\hat{\beta}}\right) 
		\notag
		\\&&
		\sim
		\frac{u^p}{n\left(ac_{\widehat{Y}}\right)^{-1/\kappa}}
		\Psi\left(u\right) 
		\notag
		\\&&\quad
		\times
		\sum_{k \in  K_{u,n}} \MB_{B_{\kappa}}\left(L\left(ac_{\widehat{Y}}\right)^{1/\kappa} \tau_{u,n,k}^{1-\kappa/\alpha},[0,n]\right)
		\notag
		\\
		&&
		\quad\quad\quad\quad\quad
		\times \exp\left(-du^{1+\left(p-2/\kappa\right)\hat{\gamma}} \tau_{u,n,k}^{\hat{\gamma}}
		-bu^{2+\left(p-2/\kappa\right)\hat{\beta}} \tau_{u,n,k}^{\hat{\beta}}\right)\left(\tau_{u,n,k+1}-\tau_{u,n,k}\right)
		\notag
		\\&&
		\sim
		\frac{u^p}{n\left(ac_{\widehat{Y}}\right)^{-1/\kappa}}
		\Psi\left(u\right)
		\notag
		\\&&\quad
		\times
		\int\limits_{\frac{1}{M'}}^{M'} 
		\exp\left(-du^{1+\left(p-2/\kappa\right)\hat{\gamma}} z^{\hat{\gamma}}
		-bu^{2+\left(p-2/\kappa\right)\hat{\beta}} z^{\hat{\beta}}\right)
		\notag
		\\
		&&
		\quad\quad\quad
		\times
		\MB_{B_{\kappa}}\left(L\left(ac_{\widehat{Y}}\right)^{1/\kappa} z^{1-\kappa/\alpha},[0,n]\right) dz.
		\EQN
		Therefore, \eqref{constant} is satisfied with 
	\BQNY
		\bar{F}(x) = \frac{
				\int\limits_{\frac{1}{M'}}^{M'}
				\exp\left(-d\mathbb{I}
				\left(1+\left(p-2/\kappa\right)\hat{\gamma}=0\right) z^{\hat{\gamma}}
				-b\mathbb{I}\left(2+\left(p-2/\kappa\right)\hat{\beta}=0\right) z^{\hat{\beta}}\right)
				A(z)
				dz
		}{
		\int\limits_{\frac{1}{M'}}^{M'}
		\exp\left(-d\mathbb{I}
		\left(1+\left(p-2/\kappa\right)\hat{\gamma}=0\right) z^{\hat{\gamma}}
		-b\mathbb{I}\left(2+\left(p-2/\kappa\right)\hat{\beta}=0\right) z^{\hat{\beta}}\right)
		A(0)
		dz
	},
	\EQNY
	where $A(z)=\MB_{B_{\kappa}}\left(x\left(ac_{\widehat{Y}}\right)^{1/\kappa} z^{1-\kappa/\alpha}\right)$.
		
		Indeed, the existence, positivity and finiteness of $\MB_{B_{\kappa}}(x)$ was proved in \cite{debicki2023sojourn}. In particular, $$\sup_{n > 1} \frac{\MB_{B_{\kappa}}\left(0,[0,n]\right)}{n}
		<\infty;$$ hence, we can swap the integral and the limit in $n$ in \eqref{supsinglesum} and \eqref{sojsinglesum}, since $$\sup_{n > 1, x \geq 0} \frac{\MB_{B_{\kappa}}\left(x,[0,n]\right)}{n}\leq 
		\sup_{n > 1} \frac{\MB_{B_{\kappa}}\left(0,[0,n]\right)}{n}
		<\infty.$$
		
		Therefore, condition A2 is also satisfied.
		
		The condition \eqref{singlesumbig} now follows from
		\eqref{supsinglesum}, so, the proof of A3 is finished, hence, \eqref{supsinglesum} implies
		\BQN
		&&
		\label{asympun}
		\pk*{\underset{t \in E\left(u,n\right)}{\sup} \frac{\widehat{X}\left(t\right)}{1+
				\frac{dt^{\hat{\gamma}}}{u}}>u}
		\notag
		\\&&\sim
		\frac{u^p}{n\left(ac_{\widehat{Y}}\right)^{-1/\kappa}}
		\Psi\left(u\right) \MB_{B_{\kappa}}\left(0,[0,n]\right)\notag
		\\&&
		\quad
		\times \int\limits_{\frac{1}{M'}}^{M'} 
		\exp\left(-du^{1+\left(p-2/\kappa\right)\hat{\gamma}} z^{\hat{\gamma}}
		-bu^{2+\left(p-2/\kappa\right)\hat{\beta}} z^{\hat{\beta}}\right) dz,
		\EQN
		which implies \eqref{A1finish}, so, the proof of A1 is also finished, and \eqref{asympun}  implies
		\BQNY
		&&\pk*{\underset{t \in E_u}{\sup} \frac{\widehat{X}\left(t\right)}{1+
				\frac{dt^{\hat{\gamma}}}{u}}>u}
		\\
		&&
		\sim\frac{u^p}{\left(ac_{\widehat{Y}}\right)^{-1/\kappa}}
		\Psi\left(u\right) \MB_{B_{\kappa}}\left(0\right)
		\int\limits_{\frac{1}{M'}}^{M'} 
		\exp\left(-du^{1+\left(p-2/\kappa\right)\hat{\gamma}} z^{\hat{\gamma}}
		-bu^{2+\left(p-2/\kappa\right)\hat{\beta}} z^{\hat{\beta}}\right) dz
		\\&&\sim\frac{u^p}{\left(ac_{\widehat{Y}}\right)^{-1/\kappa}}
		\Psi\left(u\right) \MB_{B_{\kappa}}\left(0\right)
		\\
		&&
		\quad
		\times
		\int\limits_{\frac{1}{M'}}^{M'} 
		\exp\left(-d\mathbb{I}
		\left(1+\left(p-2/\kappa\right)\hat{\gamma}=0\right) z^{\hat{\gamma}}
		-b\mathbb{I}
		\left(2+\left(p-2/\kappa\right)\hat{\beta}=0\right) z^{\hat{\beta}}\right) dz
		\EQNY

		On combining the above with \Cref{th1} and \eqref{supsinglesum}, we obtain 
		\BQNY
		&&\pk*{\int\limits_{\frac{1}{M'}u^{p-2/\kappa}}^{M'u^{p-2/\kappa}} \mathbb{I}\left(\overline{X}(t^{\kappa/\alpha})>\frac{u\left(1+\frac{dt^{\hat{\gamma}}}{u}\right)}{\sigma_{\widehat{X}}\left(t\right)}\right) t^{\kappa/\alpha-1}dt >Lu^{-\frac{2+p\left(\alpha-\kappa\right)}{\alpha}}}
		\\&&
		\sim
		\frac{u^p \Psi\left(u\right)}{\left(ac_{\widehat{Y}}\right)^{-1/\kappa}}
		\\&&\quad
		\times
		\int\limits_{\frac{1}{M'}}^{M'} 
		\exp\left(-d\mathbb{I}
		\left(1+\left(p-2/\kappa\right)\hat{\gamma}=0\right) z^{\hat{\gamma}}
		-b\mathbb{I}\left(2+\left(p-2/\kappa\right)\hat{\beta}=0\right) z^{\hat{\beta}}\right)
		\\
		&&
		\quad\quad\quad
		\times
		\MB_{B_{\kappa}}\left(L\left(ac_{\widehat{Y}}\right)^{1/\kappa} z^{1-\kappa/\alpha}\right) dz
		\EQNY
		for all $L \geq 0$. Fix $\epsilon \in [0,\frac{L}{2}]$, then
		\begin{align*}
		&\pk*{\int\limits_{0}^{\delta} \mathbb{I}\left(\overline{X}(t^{\kappa/\alpha})>\frac{u\left(1+\frac{dt^{\hat{\gamma}}}{u}\right)}{\sigma_{\widehat{X}}\left(t\right)}\right) t^{\kappa/\alpha-1}dt >L_u}
		\\ &\leq 
		\pk*{\int\limits_{\frac{1}{M'}u^{p-2/\kappa}}^{M'u^{p-2/\kappa}} \mathbb{I}\left(\overline{X}(t^{\kappa/\alpha})>\frac{u\left(1+\frac{dt^{\hat{\gamma}}}{u}\right)}{\sigma_{\widehat{X}}\left(t\right)}\right) t^{\kappa/\alpha-1}dt >\left(L-2\epsilon\right)u^{-\frac{2+p\left(\alpha-\kappa\right)}{\alpha}}}
		\\&\quad
		+\pk*{\int\limits_{0}^{\frac{1}{M'}u^{p-2/\kappa}} \mathbb{I}\left(\overline{X}(t^{\kappa/\alpha})>\frac{u\left(1+\frac{dt^{\hat{\gamma}}}{u}\right)}{\sigma_{\widehat{X}}\left(t\right)}\right) t^{\kappa/\alpha-1}dt >\epsilon u^{-\frac{2+p\left(\alpha-\kappa\right)}{\alpha}}}
		\\&\quad+
		\pk*{\int\limits_{M'u^{p-2/\kappa}}^{\delta} \mathbb{I}\left(\overline{X}(t^{\kappa/\alpha})>\frac{u\left(1+\frac{dt^{\hat{\gamma}}}{u}\right)}{\sigma_{\widehat{X}}\left(t\right)}\right) t^{\kappa/\alpha-1}dt >\epsilon u^{-\frac{2+p\left(\alpha-\kappa\right)}{\alpha}}}.
		\end{align*}
		Since for now we consider only the cases when $p>0$, by Lemmas \ref{betacut}-\ref{smallcut}
		for all $M'$ we obtain
		\BQNY
		&&\pk*{\int\limits_{0}^{\delta} \mathbb{I}\left(\overline{X}(t^{\kappa/\alpha})>\frac{u\left(1+\frac{dt^{\hat{\gamma}}}{u}\right)}{\sigma_{\widehat{X}}\left(t\right)}\right) t^{\kappa/\alpha-1}dt >L_u}
		\\&& \leq 
		\left(1+c\left(M'\right)+o\left(1\right)\right)
		\frac{u^p \Psi\left(u\right)}{\left(ac_{\widehat{Y}}\right)^{-1/\kappa}}\\
		&& 
		\quad
		\times
		\int\limits_{\frac{1}{M'}}^{M'} 
		\exp\left(-d\mathbb{I}
		\left(1+\left(p-2/\kappa\right)\hat{\gamma}=0\right) z^{\hat{\gamma}}
		-b\mathbb{I}\left(2+\left(p-2/\kappa\right)\hat{\beta}=0\right) z^{\hat{\beta}}\right)
		\\
		&&\quad\quad\quad
		\times
		\MB_{B_{\kappa}}\left(\left(L-2\epsilon\right)\left(ac_{\widehat{Y}}\right)^{1/\kappa} z^{1-\kappa/\alpha}\right) dz,
		\EQNY
		where $c\left(M'\right) \to 0$ as $M' \to \infty$.
		
		Therefore, letting $M' \to \infty$ and then letting $\epsilon \downarrow 0$ by the continuity of $\MB_{B_{\kappa}}\left(x\right)$ in $x$ and by the monotone convergence theorem we obtain
		\BQNY
		&&\pk*{\int\limits_{0}^{\delta} \mathbb{I}\left(\overline{X}(t^{\kappa/\alpha})>\frac{u\left(1+\frac{dt^{\hat{\gamma}}}{u}\right)}{\sigma_{\widehat{X}}\left(t\right)}\right) t^{\kappa/\alpha-1}dt >L_u}
		\\&& \leq 
		\left(1+o\left(1\right)\right) \frac{u^p\Psi\left(u\right)}{\left(ac_{\widehat{Y}}\right)^{-1/\kappa}}\\&&\quad
		\times
		\int\limits_{0}^{\infty} 
		\exp\left(-d\mathbb{I}
		\left(1+\left(p-2/\kappa\right)\hat{\gamma}=0\right) z^{\hat{\gamma}}
		-b\mathbb{I}\left(2+\left(p-2/\kappa\right)\hat{\beta}=0\right) z^{\hat{\beta}}\right)
		\\
		&&\quad\quad\quad
		\times
		\MB_{B_{\kappa}}\left(L\left(ac_{\widehat{Y}}\right)^{1/\kappa} z^{1-\kappa/\alpha}\right) dz.
		\EQNY
		On the other hand,
		for all $M'$
		\BQNY
		&&\pk*{\int\limits_{0}^{\delta} \mathbb{I}\left(\overline{X}(t^{\kappa/\alpha})>\frac{u\left(1+\frac{dt^{\hat{\gamma}}}{u}\right)}{\sigma_{\widehat{X}}\left(t\right)}\right) t^{\kappa/\alpha-1}dt >L_u}\\
		&&\geq 
		\pk*{\int\limits_{\frac{1}{M'}u^{p-2/\kappa}}^{M'u^{p-2/\kappa}} \mathbb{I}\left(\overline{X}(t^{\kappa/\alpha})>\frac{u\left(1+\frac{dt^{\hat{\gamma}}}{u}\right)}{\sigma_{\widehat{X}}\left(t\right)}\right) t^{\kappa/\alpha-1}dt >L_u}
		\\&& \geq 
		\left(1+o\left(1\right)\right)
		\frac{u^p \Psi\left(u\right)}{\left(ac_{\widehat{Y}}\right)^{-1/\kappa}}\\&& 
		\quad
		\times
		\int\limits_{\frac{1}{M'}}^{M'} 
		\exp\left(-d\mathbb{I}
		\left(1+\left(p-2/\kappa\right)\hat{\gamma}=0\right) z^{\hat{\gamma}}
		-b\mathbb{I}\left(2+\left(p-2/\kappa\right)\hat{\beta}=0\right) z^{\hat{\beta}}\right)
		\\
		&&\quad\quad\quad
		\times
		\MB_{B_{\kappa}}\left(L\left(ac_{\widehat{Y}}\right)^{1/\kappa} z^{1-\kappa/\alpha}\right) dz,
		\EQNY
		so, finally, if $p>0$, then (the first asymptotics follows from \eqref{deltacut})
		\BQNY
		&&\pk*{\int\limits_{0}^{T} \mathbb{I}\left(X\left(t\right)-dt^{\gamma}>u\right) dt >\left(\frac{\kappa}{\alpha}\right)L_u}\\&&\sim\pk*{\int\limits_{0}^{\delta} \mathbb{I}\left(\overline{X}(t^{\kappa/\alpha})>\frac{u\left(1+\frac{dt^{\hat{\gamma}}}{u}\right)}{\sigma_{\widehat{X}}\left(t\right)}\right) t^{\kappa/\alpha-1}dt >L_u}
		\\&&\sim \frac{u^p\Psi\left(u\right)}{\left(ac_{\widehat{Y}}\right)^{-1/\kappa}}
		\\&&\quad
		\times
		\int\limits_{0}^{\infty} 
		\exp\left(-d\mathbb{I}
		\left(1+\left(p-2/\kappa\right)\hat{\gamma}=0\right) z^{\hat{\gamma}}
		-b\mathbb{I}\left(2+\left(p-2/\kappa\right)\hat{\beta}=0\right) z^{\hat{\beta}}\right)
		\\
		&&
		\quad\quad\quad
		\times
		\MB_{B_{\kappa}}\left(L\left(ac_{\widehat{Y}}\right)^{1/\kappa} z^{1-\kappa/\alpha}\right) dz,
		\EQNY
		which, combined with the fact that $c_{\widehat{Y}}=c_Y(\kappa/\alpha)^{\kappa}$, finishes the proof for $p>0$.
		
		In the cases when $p\left(\alpha,\beta,\gamma,\kappa\right)=0$
		(that is, case (iii)), again, fix $M'>0$. 
		Let 
		$q = \min\left(-2/\kappa,-2/\hat{\beta},-1/\hat{\gamma}\right)$, then if $p=0$, then $L_u=Lu^{q\left(\kappa/\alpha\right)}$. Remark that
		
		\BQNY
		&&
		\pk*{\int\limits_{0}^{M'u^q} \mathbb{I}\left(\frac{\widehat{X}\left(t\right)}{1+\frac{dt^{\hat{\gamma}}}{u}}>u\right) t^{\kappa/\alpha-1}dt >L_u}
		\\
		&&=
		\pk*{
			\int\limits_{0}^{M'} \mathbb{I}\left(\frac{\widehat{X}\left(tu^q\right)}{1+\frac{dt^{\hat{\gamma}}}{u^{1-q\hat{\gamma}}}}>u\right) t^{\kappa/\alpha-1}dt>L_u u^{-q\left(\kappa/\alpha\right)}
		}
		\\&&=\pk*{
			\int\limits_{0}^{M'} \mathbb{I}\left(\frac{\overline{\widehat{X}}\left(tu^q\right)\sigma_{\widehat{X}}\left(tu^q\right)}{1+\frac{dt^{\hat{\gamma}}}{u^{1-q\hat{\gamma}}}}>u\right) t^{\kappa/\alpha-1}dt>L
		},
		\EQNY
		so, we can apply Lemma \ref{the-weak-conv} to the process $\frac{\overline{\widehat{X}}\left(tu^q\right)\sigma_{\widehat{X}}\left(tu^q\right)}{1+\frac{dt^{\hat{\gamma}}}{u^{1-q\hat{\gamma}}}}$ with $\eta\left(dt\right)=t^{\kappa/\alpha-1}dt$, $g_{u,j}=u$.
		
		Condition \Cref{C0} will hold obviously. 
		
		Let us check condition \Cref{C1}: 
		\BQNY
		g_{u,j}^2\left(1-\frac{\sigma_{\widehat{X}}\left(tu^q\right)}{1+\frac{dt^{\hat{\gamma}}}{u^{1-q\hat{\gamma}}}}\right)\sim u^2 b\left(tu^q\right)^{\hat{\beta}}+u d\left(tu^q\right)^{\hat{\gamma}} \to bt^{\hat{\beta}}\mathbb{I}\left(q=-2/\hat{\beta}\right)+
		dt^{\hat{\gamma}}\mathbb{I}\left(q=-1/\hat{\gamma}\right)
		\EQNY 
		as $u \to \infty$.
		
		As for \Cref{C2}:
		\BQNY
		g_{u,j}^2 \mathrm{Var}\left(
		\overline{\widehat{X}}\left(tu^q\right)
		-\overline{\widehat{X}}\left(su^q\right)
		\right)
		\sim 2aV_{\widehat{Y}}\left(t,s\right)u^{2+q \kappa} \to 2aV_{\widehat{Y}}\left(t,s\right)
		\mathbb{I}\left(q=-2/\kappa\right)
		\EQNY 
		as $u \to \infty$.
		
		Condition \Cref{C3}, as before, follows from \eqref{bddcorr}.
		
		Therefore, for $\hat{h}\left(t\right)=bt^{\hat{\beta}}\mathbb{I}\left(q=-2/\hat{\beta}\right)+
		dt^{\hat{\gamma}}\mathbb{I}\left(q=-1/\hat{\gamma}\right)$ and for a Gaussian centered process  $\zeta\left(t\right)=\widehat{Y}\left(a^{1/\kappa}t\right)\mathbb{I}\left(q=-2/\kappa\right)$ on $[0,\infty)$ we obtain
		\BQN
		\label{zerasymp}
		&&
		\pk*{\int\limits_{0}^{M'u^q} \mathbb{I}\left(\frac{\overline{\widehat{X}}\left(t\right)}{1+\frac{dt^{\hat{\gamma}}}{u}}>u\right) t^{\kappa/\alpha-1}dt >L_u}
		\notag
		\\
		&&\sim\Psi\left(u\right)\MB_{\zeta}^{\hat{h},\eta}\left(L,[0,M']\right)
		\notag
		\\
		&&=
		\Psi(u)
		\MB_{\widehat{Y} \mathbb{I}\left(q=-2/\kappa\right)}^{\hat{h}(a^{-1/\kappa} \cdot),\eta}\left(L a^{1/\kappa},[0,M' a^{1/\kappa}]\right)
		\EQN
		as $u \to \infty$.
		
		By \eqref{deltacut} and Lemmas \ref{betacut}-\ref{kappacut}, 
		letting $M' \to \infty$ in the above asymptotics, we obtain
		\BQNY
		&&
		\pk*{\int\limits_{0}^{T} \mathbb{I}\left(X\left(t\right)-dt^{\gamma}>u\right) dt >\left(\frac{\kappa}{\alpha}\right)L_u}
		\\
		&&\sim\pk*{\int\limits_{0}^{\delta} \mathbb{I}\left(\overline{X}(t^{\kappa/\alpha})>\frac{u\left(1+\frac{dt^{\hat{\gamma}}}{u}\right)}{\sigma_{\widehat{X}}\left(t\right)}\right) t^{\kappa/\alpha-1}dt >L_u}
		\\&&\sim
		\Psi\left(u\right)\MB_{\widehat{Y} \mathbb{I}\left(q=-2/\kappa\right)}^{\hat{h}(a^{-1/\kappa}\cdot),\eta}
		\left(L a^{1/\kappa},[0,\infty)\right)
		\\
		&&=\Psi\left(u\right)\MB_{Y\mathbb{I}\left(q=-2/\kappa\right)}^{h}
		\left(L\left(\frac{\kappa}{\alpha}\right)a^{1/\kappa},[0,\infty)\right),
		\EQNY
		with $h\left(t\right)=a^{-\beta/\alpha}bt^{\beta}\mathbb{I}\left(q=-2/\hat{\beta}\right)+
		a^{-\gamma/\alpha}
		dt^{\gamma}\mathbb{I}\left(q=-1/\hat{\gamma}\right)$, and Theorem \ref{mainth} is proven.
		We do not need to check the finiteness of this Berman constant, since 
		\BQNY
		\pk*{\int\limits_{0}^{\delta} \mathbb{I}\left(\overline{X}(t^{\kappa/\alpha})>\frac{u\left(1+\frac{dt^{\hat{\gamma}}}{u}\right)}{\sigma_{\widehat{X}}\left(t\right)}\right) t^{\kappa/\alpha-1}dt >L_u}=O\left(\Psi\left(u\right)\right)
		\EQNY
		by Lemmas \ref{betacut}-\ref{kappacut} and \eqref{zerasymp}. 
	\end{proof}

	\section*{Acknowledgements}
	Financial support by SNSF Grant 200021–196888 is kindly acknowledged. 
	
	\section*{Disclosure statement}
	The author has no conflicts of interest to declare that are relevant to the content of this article.
	
	\section*{Funding}
	This work was supported by the SNSF under Grant 200021–196888. 

	\bibliographystyle{tfs}
	\bibliography{EEEA}

	\end{document}